\documentclass[a4paper]{amsart}
\usepackage{a4wide}
\usepackage{graphicx}
\usepackage{amsmath}
\usepackage{amssymb}
\usepackage{amsfonts}
\usepackage{amsthm}
\usepackage{eucal}
\usepackage{tikz}
\usepackage[colorlinks=true, linkcolor={blue!50!black}, pdfhighlight=/O, 
ocgcolorlinks=true]{hyperref}
\usepackage[normalem]{ulem}
\usepackage{comment}

\usetikzlibrary{trees, positioning, calc, cd}

\theoremstyle{plain}
\newtheorem{theorem}{Theorem}[section]

\newtheorem{lemma}[theorem]{Lemma}
\newtheorem{proposition}[theorem]{Proposition}

\newtheorem{corollary}[theorem]{Corollary}

\newtheorem*{proposition*}{Proposition}
\newtheorem*{corollary*}{Corollary}
\newtheorem*{theorem*}{Theorem}
\newtheorem{observation}[theorem]{Observation}

\theoremstyle{definition}
\newtheorem{definition}[theorem]{Definition}

\newtheorem{remark}[theorem]{Remark}

\makeatletter
\@addtoreset{theorem}{section}
\@addtoreset{subsection}{section}
\makeatother

\newcommand{\Hom}{\mathrm{Hom}}    
\newcommand{\sSet}{{\sfix\Set}}      
\newcommand{\sfix}{\mathrm{s}}	   
\newcommand{\Set}{\mathrm{Set}}    
\newcommand{\mSet}{{\sfix\Set^+}}      
\newcommand{\F}{\Hom^{\simp}}        
\newcommand{\FHom}{\Hom}	
\newcommand{\unmarked}{\flat}	   
\newcommand{\marked}{\#}		   
\newcommand{\Fun}{\mathrm{Fun}}    
\newcommand{\Cat}{\mathrm{Cat}}    
\newcommand{\cnerv}{{\mathrm N^\mathrm{c}}} 
\newcommand{\hocolim}{\mathrm{hocolim}} 

\newcommand{\C}{\mathcal C}     
\newcommand{\CC}{\mathfrak C}     
\newcommand{\DD}{\mathfrak D}     
\newcommand{\M}{\mathfrak M}      
\newcommand{\simp}{\mathrm{s}}    
\newcommand{\Catsimp}{\Cat^{\sSet}} 
\newcommand{\nerve}{\mathrm{N}}     
\newcommand{\NL}{\mathrm {NL}}     
\newcommand{\Nec}{\mathrm{Nec}} 
\DeclareMathOperator{\colim}{colim} 
\newcommand{\op}{\mathrm{op}}
\newcommand{\Strun}{\mathrm{Str}}
\newcommand{\Unun}{\mathrm{Un}}
\newcommand{\Str}{{\Strun^{+}}}
\newcommand{\Un}{{\Unun^+}}
\newcommand{\Nat}{\mathrm{Nat}}

\usepackage{pict2e,picture}     
\makeatletter
\DeclareRobustCommand{\bbDelta}{{\mathpalette\bb@Delta\relax}}
\newcommand{\bb@Delta}[2]{%
  \begingroup
  \sbox\z@{$\m@th#1\Delta$}%
  \dimendef\Dht=6 \dimendef\Dwd=8
  \setlength{\Dwd}{\wd\z@}%
  \setlength{\Dht}{\ht\z@}%
  \begin{picture}(\Dwd,\Dht)
  \put(0,0){$\m@th#1\Delta$}
  \put(.52\Dwd,.77\Dht){\line(-13,-26){.35\Dht}}
  \end{picture}%
  \endgroup
}
\makeatother

\newcommand*\cocolon{%
        \nobreak
        \mskip6mu plus1mu
        \mathpunct{}%
        \nonscript
        \mkern-\thinmuskip
        {:}%
        \mskip2mu
        \relax
}

\begin{document}

\title{A short proof of the straightening theorem}

\author{Fabian Hebestreit}
\address{Mathematisches Institut, WWU Münster, Germany}
\email{f.hebestreit@uni-muenster.de}

\author{Gijs Heuts}
\address{Mathematisch Instituut, Universiteit Utrecht, The Netherlands}
\email{g.s.k.s.heuts@uu.nl}

\author{Jaco Ruit}
\address{Mathematisch Instituut, Universiteit Utrecht, The Netherlands}
\email{j.c.ruit@uu.nl}
\date{}

\begin{abstract}
We provide a short and reasonably self-contained proof of Lurie's straightening equivalence, relating cartesian fibrations over a given $\infty$-category $S$ with contravariant functors from $S$ to the $\infty$-category of small $\infty$-categories.
\end{abstract}

\maketitle

\setcounter{tocdepth}{1}
\tableofcontents

\section{Introduction}

Lurie's straightening theorem provides an equivalence
\[\begin{tikzcd} \Strun \colon \mathrm{Cart}(S) \ar[shift left]{r} & \ar[shift left]{l} \mathrm{Fun}(S^\mathrm{op},\mathrm{Cat}_\infty) \cocolon \Unun\end{tikzcd}\] 
for any $\infty$-category $S$, where the left hand side denotes the $\infty$-category of cartesian fibrations over $S$ and $\mathrm{Cat}_\infty$ is the $\infty$-category of small $\infty$-categories \cite{HTT}. Along with Joyal's proof of Grothendieck's homotopy hypothesis \cite{Joyalqcat}, it is one of the cornerstones of higher category theory. Cartesian fibrations are the $\infty$-categorical generalisation of fibred categories in ordinary category theory and the straightening equivalence generalises Grothendieck's equivalence between fibred categories and prestacks, as well as the correspondence between covering spaces and sets with an action of the fundamental group of the base, and also the classification of Kan or Serre fibrations in homotopy theory by maps into a classifying space. Its applications range from the $\infty$-categorical version of Yoneda's lemma to the computation of colimits in the $\infty$-categories of spaces and of small $\infty$-categories; it is also baked into the foundations of Lurie's treatment of symmetric monoidal $\infty$-categories and higher operads. In fact, the coherence issues involved in constructing functors of $\infty$-categories into $\mathrm{Cat}_\infty$ (or the $\infty$-category of spaces/$\infty$-groupoids) are in practice almost always solved by writing down the corresponding (co)cartesian fibrations instead.

It is therefore somewhat surprising that the only complete proof of the theorem in the literature at the time of writing is the original one in \cite{HTT}.\footnote{During the revision process of this manuscript Cisinski and Nguyen in fact gave another full proof in \cite{CisNgucurrent}.} Lurie's approach starts by writing both sides as $\infty$-categories associated to simplicial model categories via the coherent nerve $\cnerv$ (applied to the bifibrant objects); recall that this is a functor $\cnerv \colon \Cat^{\sSet} \rightarrow \sSet$ that assigns an $\infty$-category to every category enriched in Kan complexes. The left-hand side of the straightening equivalence underlies the cartesian model structure on the category $\mSet/S^\sharp$ of marked simplicial sets over $S$, and the right hand side underlies the projective model structure on the category of simplicially enriched functors $\CC(S)^\op \rightarrow \mSet$, where $\CC$, the path-category functor, is the left adjoint to the coherent nerve and $\mSet$ is equipped with the marked Joyal model structure; it satisfies $\cnerv(\mSet) \simeq \mathrm{Cat}_\infty$. For the left hand side this translation is essentially a triviality and for the right hand side it is an instance of a general principle: Lurie shows that for any sufficiently nice simplicial model category $\mathfrak M$ the category of simplically enriched functors $\CC(S) \rightarrow \mathfrak M$ has $\Fun(S,\cnerv(\mathfrak M))$ as its underlying $\infty$-category.

The present paper is concerned with the remaining step, which is the most arduous in Lurie's treatment. Namely, we provide a new proof of the following result:

\begin{theorem*}
For any $\infty$-category $S$ the marked straightening-unstraightening adjunction
\[
\begin{tikzcd}
\Str \colon \mSet/S^\sharp \ar[shift left]{r} & \mathrm{Fun}^\simp(\mathfrak{C}(S)^{\mathrm{op}},\mSet) \ar[shift left]{l} \cocolon \Un
\end{tikzcd}
\]
is a Quillen equivalence for the cartesian model structure on the left and the projective model structure based on the marked Joyal model structure on the right. 
\end{theorem*}

As already observed by Lurie \cite[Subsection 3.2.4]{HTT}, the Quillen equivalence above is not itself simplicial; the straightening functor is not even simplicially enriched. However, the unstraightening functor is simplicial and this is good enough to obtain an equivalence on coherent nerves as desired.

Lurie's strategy of proof consists of two main steps. First, he treats the basic case of a simplex $S = \Delta^n$ and then establishes an induction to get from simplices $\Delta^n$ to a general simplicial set $S$. Both of these steps are nontrivial: the first involves a rather intricate combinatorial analysis using `mapping simplices', whereas the second relies on a delicate analysis of homotopy colimits of simplicially enriched categories.

By contrast, our proof works for a fixed $\infty$-category $S$. Rather than an induction on the base, we provide an inductive argument on the `total space' of a cartesian fibration $X \rightarrow S$: since cartesian equivalences between cartesian fibrations can be detected fibrewise, it follows from basic properties of the construction that $\mathbf R\Un$ detects weak equivalences between fibrant objects, so by general abstract nonsense it suffices to show that the derived unit
\[\eta \colon p \longrightarrow \mathbf R\Un\mathbf L\Str(p)\]
is a cartesian equivalence for every cartesian fibration $p \colon E \rightarrow S$ (marked by its cartesian edges). The key insight that allows an inductive proof of this statement is that the right derived functor $\mathbf{R}\Un$ preserves homotopy colimits. Our proof of this result relies on several facts about the cartesian model structure (in particular, the behavior of this model structure under `base change' along a map $T \rightarrow S$) that Lurie originally deduced from the straightening theorem, but which have since been given independent treatments in the literature, primarily by Cisinski and Nguyen \cite{Cisinskibook, Nguyenthesis, Nguyen}; for the reader's convenience we include an account of the requisite results. Let us also immediately mention that our argument is geared to work for $S$ an $\infty$-category (as reflected in the statement of the theorem above) but for completeness' sake we explain how to deduce the case of a general base from this. 

The base case of our induction is the case of a map $\Delta^0 = E \xrightarrow{p} S$. A fibrant replacement for $p$ is given by the slice projection $S_{/y} \rightarrow S$, where $y \in S_0$ is the value of $p$. Unwinding definitions now translates the assertion on $\eta$ in this case to the natural map
\[\Hom_S(x,y) \longrightarrow \Hom_{\mathbf R\cnerv\mathbf L\CC(S)}(x,y)\]
being a homotopy equivalence for all $x \in S$. This in turn is a direct consequence of the fact that the functors 
\[
\begin{tikzcd}
\CC \colon \sSet \ar[shift left]{r} & \Cat^\sSet \ar[shift left]{l} \cocolon \cnerv
\end{tikzcd}
\]
form a Quillen equivalence between Joyal's model structure on the left and Bergner's model structure on the right. Lurie deduces this result from the straightening theorem (in the case of right fibrations) using the line of thought above in reverse, and another proof, comparing both sides to complete Segal spaces, was given by Joyal \cite{Joyalsimpcats}. A much more direct argument for the equivalence of these mapping spaces was given by Dugger and Spivak in \cite{duggerspivakmapping,duggerspivakrigid} and for the reader's convenience we again give a brief account of the part of their work that is relevant for us.

Let us conclude the introduction by briefly commenting on the relation to the literature. The second author and Moerdijk gave a simplified proof of straightening of left/right fibrations and functors into $\mathrm{Gpd}_\infty \subset \Cat_\infty$ in \cite{left1,left2} by first establishing an equivalence $\mathbf R\Unun \simeq \mathbf Lh$ for $h \colon \Fun^\simp(\CC(S)^\op,\sSet) \rightarrow \sSet/S$ a variant of the usual homotopy colimit functor. Observing that this implies that $\mathbf R\Unun$ preserves homotopy colimits, the remainder of their argument is (at least in spirit) not so different from the inductive one we give here. The present paper in fact grew out of an attempt to generalise this strategy to the case of (co)cartesian fibrations. Ultimately we realised that one can prove directly that $\mathbf R\Un$ preserves homotopy colimits and this is the key innovation of the present work. Specialised to left fibrations one could thus regard our approach here as a further simplification of that in \cite{left2}.

Straightening of left fibrations was also treated in detail by Cisinski in \cite{Cisinskibook}. There is, however, a crucial difference in that Cisinksi's arguments do not pass through simplicially enriched categories, but are rather based on a different definition of the subcategory $\mathrm{Gpd}_\infty \subset \Cat_\infty$, which essentially trivialises the statement at the point-set level. While very elegant, this argument outsources a significant amount of work to establishing the relation between the categories of $\infty$-groupoids and Kan complexes, which is a simple application of Joyal's lifting theorem \cite[Proposition 1.2.4.3]{HTT} in Lurie's approach. Finally in \cite{Kerodon}, mimicking parts of Cisinski's approach but keeping the definition of $\Cat_\infty$ in terms of the coherent nerve, Lurie gave a new proof of the straightening equivalence at the level of homotopy categories. He also promised to get back to the full statement in the future and we hope our proof will make the wait more bearable.

\subsection*{Organisation}
In Section \ref{sec:preliminaries} we collect the input that we require for the proof of the straightening theorem. The proof itself is then the content of the short Section \ref{sec:proof}. Largely for the convenience of the reader, Section \ref{sec:modelstructures} contains a basic discussion of the cartesian and projective model structure. Section \ref{sec:basechange} and \ref{sec:Necklaces} finally take up the task of verifying the various inputs for our proof as stated in the second section. As mentioned in the introduction, these are mostly well-known facts but many were originally derived from the straightening theorem in \cite{HTT,HA}. To remove the risk of circularity we provide detailed a priori verifications, substantial parts of which we extract from the existing literature: Section \ref{sec:basechange}, concerning base change results, closely follows work of Cisinski and Nguyen \cite{Cisinskibook, Nguyen}, and Section \ref{sec:Necklaces} concerning path categories is a stripped-down account of work of Dugger and Spivak \cite{duggerspivakrigid}. We include this material as a service to the reader and to provide an efficient path to obtaining these results, that one could for example use in a lecture course (as we have).

\subsection*{Acknowledgments}
It is a pleasure to acknowledge fruitful exchanges with Bastiaan Cnossen, Markus Land, Sil Linskens, Ieke Moerdijk and Ferdinand Wagner about the contents of this paper. Furthermore, we would like to thank all participants of the lecture course `Higher categories and homotopical algebra' held in the summer term 2020 at the University of Bonn for sitting through the first exposition of the material. 

During the preparation of this paper FH was a member of the Hausdorff Center for Mathematics at the University of Bonn funded by the German Research Foundation (DFG), grant no.\ EXC 2047-390685813, GH was supported by the European Research council (ERC) through the grant ``Chromatic homotopy theory of spaces'', grant no. 950048, and JR was funded by the Dutch Research Council (NWO) through the grant ``The interplay of orientations and symmetry'', grant no. OCENW.KLEIN.364. During the revision FH was further supported by the DFG through the CRC/TRR 358 ``Integral structures in representation theory and geometry", grant no. TRR 358-49139240, at the University of Bielefeld.

\section{Preliminaries and recollections}
\label{sec:preliminaries}

In this section we recall the basic definitions of the objects involved in the straightening theorem. At the same time we state a number of results that we shall use in its proof. Most of these are contained in \cite{HTT}, but as mentioned are derived from the straightening equivalence there and have since been given independent treatments in the literature. To spare the reader a treasure hunt and demonstrate that our arguments are non-circular, we discuss the proofs of all of these results, either giving a self-contained argument or a precise reference that contains a proof independent of the straightening equivalence.

\subsection{Path categories and the coherent nerve}\label{sec:pathcatnerv}
\label{subsec:coherentnerve}

We will briefly recall the coherent nerve and its adjoint, taking a simplicial set to what we shall call its path category (following \cite{Kerodon}, \cite[Definition 1.1.5.5]{htt}). We start by setting up some notation. Given a cocomplete category $\C$ and a cosimplicial object $X \colon \bbDelta \rightarrow \C$, there is an induced adjunction
\[
\begin{tikzcd}
{|\cdot|_X} \colon \sSet \ar[shift left]{r} & \C \ar[shift left]{l} \cocolon {\mathrm{Sing}_X} 
\end{tikzcd}
\]
where $\mathrm{Sing}_X(c) = [n \mapsto \Hom_\C(X_n,c)]$ and its left adjoint $|\cdot|_X$ is the left Kan extension of $X$ along the Yoneda embedding $\Delta \colon \bbDelta \rightarrow \sSet$.

The coherent nerve arises in this way for a certain cosimplicial object $\CC \colon \bbDelta \rightarrow \Catsimp$ in the category of simplicially enriched categories. To construct it, let $I$ be a finite totally ordered set $I$ with $i \leq j \in I$ and consider the poset
\[\mathrm{P}^I_{i,j} := \{A \subseteq [i,j] \mid i,j \in A\}\] 
ordered by inclusion (we shall drop the superscript $I$ when there is no danger for confusion). Then define $\CC_I \in \Catsimp$ as having object set $I$ and
\[\F_{\CC_I}(i,j) = \nerve(\mathrm{P}_{i,j}),\]
where $\nerve \colon \Cat \rightarrow \sSet$ is the nerve functor (i.e., the right adjoint arising from the inclusion $\bbDelta \subset \Cat$); to avoid confusion later we shall generally write $\F_\mathfrak D(x,y)$ for the simplicial set of morphisms in a simplicially enriched category $\mathfrak D$ (whereas $\Hom_\mathfrak D(x,y)$ will always refer to its $0$-simplices). The composition in $\CC_I$ is induced by the operation $\mathrm{P}_{j,k} \times \mathrm{P}_{i,j} \rightarrow \mathrm{P}_{i,k}$ taking unions of subsets. One defines path categories and coherent nerves
\[
\begin{tikzcd}
{\CC} \colon \sSet \ar[shift left]{r} & \Catsimp \ar[shift left]{l}\cocolon {\cnerv},
\end{tikzcd}
\]
via the same procedure with $\cnerv = \mathrm{Sing}_\CC$; here we have abbreviated the notation $|\cdot|_{\CC}$ for the left adjoint by $\CC$. 

We will need a good understanding of the relation between mapping spaces in an $\infty$-category $S$ and mapping spaces in the corresponding path category $\CC(S)$. The simplicial sets $\F_{\CC(S)}(s,t)$ admit a convenient and very explicit description in terms of \emph{necklaces} in the simplicial set $S$. This description is due to Dugger--Spivak \cite{duggerspivakmapping}; we will review it in Section \ref{sec:Necklaces} below. For mapping spaces in the $\infty$-category $S$ itself there are various standard models. The simplest arise by first picking one of the three `suspension functors' $\Sigma \colon \sSet \rightarrow \sSet$ given by
\[X \longmapsto (\Delta^0 * X)/X, \quad (X \times \Delta^1)/\{X \times 0, X \times 1\}, \quad (X * \Delta^0)/X\]
according to the cases of the left, symmetric and right mapping spaces, respectively; then the corresponding mapping space $\mathrm{Hom}_S(s,t)$ has $n$-simplices given by maps $\Sigma \Delta^n \rightarrow S$ sending the initial vertex $0$ to $s$ and the final vertex $1$ to $t$. For definiteness we will use the right mapping spaces and we will write $\mathrm{Hom}^\mathrm R_S(s,t)$ to emphasise this choice, when this becomes relevant. We will, however, also have to use the left mapping spaces in Section \ref{subsec:leftfullyfaith} and denote these by $\mathrm{Hom}^\mathrm L_S(s,t)$. 

Let us write $Q^n$ for the simplicial set $\F_{\CC(\Sigma \Delta^n)}(0,1)$. Note that for a simplicial category $\mathfrak D$ and objects $c,d$ of it, the mapping spaces of its coherent nerve are given by
\[\Hom_{\cnerv\mathfrak D}(c,d) \cong \mathrm{Sing}_Q(\F_\mathfrak D(c,d)),\]
essentially by adjunction.
This follows from the observation that the natural inclusion
$$
\Hom_{\sSet}(\Delta^n, \Hom_{\cnerv\mathfrak D}(c,d)) \rightarrow \Hom_{\sSet}(\Sigma \Delta^n, \cnerv\mathfrak D) \cong \Hom_{\Cat^\sSet}(\mathfrak{C}(\Sigma \Delta^n), \mathfrak D)
$$ 
selects those maps $\mathfrak{C}(\Sigma \Delta^n) \rightarrow \mathfrak{D}$  that carry the (only) objects $0$ and $1$ to $c$ and $d$ respectively. 
It thus follows that
the image of the above inclusion is isomorphic to
$$
\Hom_{\Cat^\sSet}(\F_{\CC(\Sigma \Delta^n)}(0,1), \F_{\mathfrak{D}}(c,d)) = \Hom_{\Cat^\sSet}(Q^n, \F_{\mathfrak{D}}(c,d)).
$$ 
Now in particular, the unit map $S \rightarrow \cnerv\CC(S)$ provides a canoncial map
\[\Hom_S(s,t) \longrightarrow \mathrm{Sing}_Q(\F_{\CC(S)}(s,t)).\]
For the particular case of right mapping spaces (corresponding to the third choice of suspension functor above) we will need the following result:

\begin{theorem}\label{thm:duggerspivak}
Let $S$ be an $\infty$-category and $s,t \in S_0$. Then the adjoint
\[|\Hom_S(s,t)|_Q \longrightarrow \F_{\CC(S)}(s,t)\]
of the map constructed above is a weak homotopy equivalence. 
\end{theorem}

As mentioned in the introduction the first proof of this result is due to Lurie, who derives it from the straightening equivalence, see \cite[Section 2.2]{HTT} and particularly \cite[Proposition 2.2.4.1]{HTT}. This dependence was eliminated by Joyal in \cite{Joyalsimpcats}, who based his argument on the equivalence between $\infty$-categories and complete Segal spaces. A much more direct proof is given in by Dugger and Spivak in \cite{duggerspivakrigid}. We review how their arguments yield Theorem \ref{thm:duggerspivak} in Section \ref{sec:Necklaces} below, see Theorem \ref{thm:duggerspivakrestate} specifically.

\subsection{The contravariant and cartesian model structures}
\label{subsec:covariantcocartmodelstruct}

In the present section we recall the definitions of the contravariant and cartesian model structures. The former is originally due to Joyal \cite{Joyalqcatappl} and the extension to the latter to Lurie \cite{HTT}. The proofs of these results are independent of the straightening theorem; particularly streamlined treatments have been given by Cisinski \cite[Section 4.1]{Cisinskibook} and Nguyen \cite{Nguyen,Nguyenthesis}. We will review these proofs in some detail in Section \ref{subsec:cocartesianmodelstruct} below, as we need some particular ingredients from them.

\begin{definition}
A map among simplicial sets is called \emph{right anodyne} if it lies in the weakly saturated class generated by the right horn inclusions $\Lambda_i^n \rightarrow \Delta^n$ with $0 < i \leq n$. A map is a \emph{right fibration} if it has the right lifting property with respect to right anodynes.
\end{definition}

There is of course a dual version of this definition; the weakly saturated class of left anodynes is generated by the left horn inclusions $\Lambda_i^n \rightarrow \Delta^n$ with $0 \leq i < n$ and the left fibrations are the maps having the right lifting property with respect to these.

\begin{theorem}[Joyal]
There exists a model structure on $\sSet/S$ characterised by the fact that its cofibrations are the injections and its fibrant objects are the right fibrations $X \rightarrow S$. This model structure has the following further properties:
\begin{itemize}
\item[(1)] All fibrations are right fibrations and dually all right anodynes are trivial cofibrations. The converse to both statements holds for maps with fibrant target (but not in general).
\item[(2)] A map 
\[
\begin{tikzcd}
X \ar{rr}{f}\ar{dr}[swap]{p} && Y \ar{dl}{q} \\
& S &
\end{tikzcd}
\]
between right fibrations $p$ and $q$ over $S$ is a weak equivalence if and only if for every vertex $s \in S_0$, the induced map of fibres $X_s \rightarrow Y_s$ is a weak homotopy equivalence of simplicial sets.
\end{itemize}
\end{theorem}

The model structure of the theorem is called the \emph{contravariant model structure}. There is a dual version, the \emph{covariant model structure}, that has right anodynes and fibrations replaced by left anodynes and fibrations, respectively. The original proof of this result is \cite[Section 8]{Joyalqcatappl} and a particularly elegant treatment is Cisinski's in \cite[Section 4.1]{Cisinskibook}, and both treatments are entirely independent of the straightening equivalence.

Let us now discuss the cartesian analogue. We use the notation
\[(-)^\unmarked \colon \sSet \longrightarrow \mSet, (-)_\unmarked \colon \mSet \longrightarrow \sSet, (-)^\marked \colon \sSet \longrightarrow \mSet, (-)_\marked \colon \mSet \longrightarrow \sSet\]
for the functors that mark only the degenerate edges, forget the markings, mark all edges, and finally extract the largest subobject all of whose edges are marked, respectively. They are adjoint to each other in the order indicated. We shall often refer to objects and morphisms in the image of $(-)^\unmarked$ as unmarked, and generally suppress degenerate edges when specifying markings. In the following definition we use the symbol $\boxtimes$ to denote the pushout-product of two maps.

\begin{definition}\label{def:markedanodyne}
A map of marked simplicial sets is called a \emph{marked injection} if its underlying map is an injection of simplicial sets. We shall call it \emph{marked right anodyne} if it lies in the weakly saturated class generated by the maps $(A^\marked \rightarrow B^\marked) \boxtimes (C \rightarrow D)$ with $A \rightarrow B$ right anodyne (in $\sSet$) and $C \rightarrow D$ a marked injection. Finally, we will call a map \emph{cartesian anodyne} if it lies in the saturation of the marked right anodynes, unmarked inner anodynes and the map $J^\unmarked \rightarrow (J,\{0 \rightarrow 1\})$, where $J$ is the nerve of the free-living isomorphism with objects $\{0,1\}$. A \emph{marked trivial fibration} is a map with the right lifting property against marked injections and a \emph{marked cartesian fibration} is a map with the right lifting property against cartesian anodynes. 
\end{definition}

Note that the map of simplicial sets underlying a cartesian anodyne is always right anodyne (since the right anodynes are closed under pushout-products with injections \cite[Corollary 2.1.2.7]{HTT}). As before there are evident dual notions of \emph{marked left anodyne} and \emph{cocartesian anodyne} morphisms, as well as \emph{marked cocartesian fibrations}.

\begin{remark}
\label{rmk:cartterminology}
Our terminology is slightly non-standard here: Lemma \ref{lem:genmarkright} shows that our marked right anodynes agree with Nguyen's \emph{cellular marked right anodynes}, and Lemma \ref{lem:charofmarkedcartfibs} shows that our marked cartesian fibrations agree with Lurie's \emph{marked fibrations}, and it follows that our cartesian anodynes agree with Lurie's \emph{marked anodynes}. Since we need to be able to consistently distinguish left and right, we have opted for the terminology introduced above.
\end{remark}

\begin{theorem}[Lurie]\label{cocartmodelstr}
Let $S$ be a marked simplicial set. There exists a model structure on $\mSet/S$ characterised by the fact that its cofibrations are the marked injections and its fibrant objects are the marked cartesian fibrations $X \rightarrow S$. This model structure has the following further properties:
\begin{itemize}
\item[(1)] All fibrations are marked cartesian fibrations and dually all cartesian anodynes are trivial cofibrations. The converse to both statements holds for maps with fibrant target (but not in general).
\item[(2)] If all edges of $S$ are marked, then a map $f\colon X \rightarrow S$ is a marked cartesian fibration precisely if the underlying map $f_\unmarked$ is a cartesian fibration of simplicial sets and the marked edges of $X$ are exactly the $f$-cartesian edges. In this case, a map 
\[
\begin{tikzcd}
X \ar{rr}{f}\ar{dr}[swap]{p} && Y \ar{dl}{q} \\
& S &
\end{tikzcd}
\]
between marked cartesian fibrations $p$ and $q$ over $S$ is a cartesian equivalence if and only if for every vertex $s \in S_0$, the induced map of fibres $(X_s)_\unmarked \rightarrow (Y_s)_\unmarked$ is a categorical equivalence of simplicial sets, i.e., a weak equivalence in the Joyal model structure.
\end{itemize}
\end{theorem}

We shall refer to the model structure of the theorem as the \emph{cartesian model structure} on $\mSet/S$. There is a dual \emph{cocartesian model structure} that has the marked cocartesian fibrations as its fibrant objects. The original reference is \cite[Section 3.1.3]{HTT}, and an elegant argument parallel to Cisinski's treatment of the unmarked case was given in \cite{Nguyen}. We will discuss some aspects of the proof in Section \ref{subsec:cocartesianmodelstruct}; in particular, we indicate how Nguyen's reliance on Smith's general existence theorem for combinatorial model structures can be avoided in the case at hand, see Proposition \ref{lem:countable} specifically. None of the arguments in either treatment rely on the straightening equivalence.

The weak equivalences in the cartesian model structure can be characterised as those maps
\[\begin{tikzcd}
E \ar[rr,"f"]\ar[rd] && E'\ar[ld] \\
& S &
\end{tikzcd}\]
such that for every marked cartesian fibration $X \rightarrow S$ the induced map 
\[\FHom^\sharp_S(E,X) \longrightarrow \FHom^\sharp_S(E',X)\]
is a homotopy equivalence; here $\FHom^\sharp_S(E,X)$ denotes the simplicial subset of $\F_{S}(E,X)$ spanned by the vertices representing marked maps (over $S$) and by the edges representing marked maps $(\Delta^1)^\sharp \times E$ (over $S$). If $X \rightarrow S$ is a marked cartesian fibration this is indeed a Kan complex by closure under pushout-products of the cartesian anodynes with marked cofibrations, see \cite[Proposition 3.1.2.3]{HTT} and the discussion in \cite[Section 3.1.3]{HTT}. Furthermore, the functor $\FHom_S^\sharp$ enriches $\mSet/S$ in simplicial sets, and indeed it is not difficult to see that $\mSet/S$ in total becomes a Kan simplicial model category as in \cite[Section 9.1]{Hirschhorn}, see either \cite[Section 3.1.4]{HTT} for a direct verification or the discussion after Proposition \ref{lem:charofmarkedcartfibs} below.

Let us also record the following observation about the `absolute' cartesian model structure, i.e., the case where $S$ is a point. It is easily proved directly or seen as a special case of Lemma \ref{lem:sliceoverqcat} below. We will refer to the model structure of the following lemma as the \emph{marked Joyal model structure} on the category $\mSet$, and to its equivalences as the \emph{marked categorical equivalences}.

\begin{lemma}
\label{lem:absolutemarked}
In case $S = \Delta^0$, the cartesian and cocartesian model structure on $\mSet/S = \mSet$ coincide. A marked simplicial set $X$ is fibrant precisely if the underlying simplicial set $X_\flat$ is an $\infty$-category and the marked edges are the equivalences in $X_\flat$. The adjoint pair
\[
\begin{tikzcd}
{(-)^\flat}  \colon \sSet \ar[shift left]{r}& \mSet \ar[shift left]{l} \cocolon {(-)_\flat}
\end{tikzcd}
\]
is a Quillen equivalence between the Joyal model structure and the marked Joyal model structure.
\end{lemma}

\subsection{Base change for (co)cartesian fibrations}\label{sec:basechangecocart}

One of our key inputs is an understanding of the behavior of the cartesian model structure under base change. A first simple observation is the following. Let $F\colon \mSet/S \rightarrow \mathcal{E}$ be a left adjoint functor that preserves cofibrations. If $F$ sends marked anodynes to weak equivalences, then it is left Quillen. Indeed, its right adjoint will preserve fibrations between fibrant objects and this suffices by \cite[Proposition E.2.14]{Joyalqcatappl}. We immediately find:

\begin{proposition}
\label{prop:basechange}
Let $f\colon T \rightarrow S$ be a map of marked simplicial sets. Then the adjoint pair
\[
\begin{tikzcd}
{f_!} \colon \mSet/T \ar[shift left]{r} & \mSet/S \ar[shift left]{l} \cocolon {f^*}
\end{tikzcd}
\]
is a Quillen adjunction with respect to the cartesian model structures.
\end{proposition}

Here $f_!$ is the functor that composes with $f$, with right adjoint $f^*$ forming the pullback along $f$. The latter functor admits a further right adjoint $f_*$. What is less obvious than Proposition \ref{prop:basechange} is that $f^*$ is also a \emph{left} Quillen functor in good cases. We shall need the following:

\begin{theorem}[Nguyen]
\label{thm:leftfibrightQuillen}
If $f\colon T \rightarrow S$ is a left fibration, then 
\[
\begin{tikzcd}
{f^*} \colon \mSet/S^\sharp \ar[shift left]{r} & \mSet/T^\sharp \ar[shift left]{l} \cocolon {f_*}
\end{tikzcd}
\]
is a Quillen adjunction with respect to the cartesian model structures.
\end{theorem}

This result originally appeared in Nguyen's thesis \cite[Proposition 3.2.6]{Nguyenthesis}: To prove it (which we do in Section \ref{sec:basechange}, see Theorem \ref{thm:leftfibrightQuillenrestate} specifically) it suffices to check that $f^*$ sends cartesian anodynes to weak equivalences. For the case of marked right anodynes this is a consequence of the following result:

\begin{proposition}\label{thm:Cartfibmarkedanod}
For $T \rightarrow S$ right anodyne and $X \rightarrow S$ a left fibration, the morphism $T \times_{S} X \rightarrow X$ in $\sSet$ is right anodyne. Similarly, for $T \rightarrow S$ marked right anodyne and $X \rightarrow  S$ a marked cocartesian fibration the morphism $T \times_{S} X \rightarrow X$ in $\mSet$ is cartesian anodyne.
\end{proposition}

The first statement is due to Joyal and the marked case is again due to Nguyen \cite[Theorem 4.45]{Nguyen}. The proof is essentially the same in both cases and is independent of the straightening equivalence; we include it in Section \ref{sec:basechange} below, see Proposition \ref{thm:Cartfibmarkedanodrestate} specifically. Note that the unmarked version already suffices to establish that
\[
\begin{tikzcd}
{f^*} \colon \sSet/S \ar[shift left]{r} & \sSet/T \ar[shift left]{l}\cocolon {f_*}
\end{tikzcd}
\] 
is a Quillen pair with respect to the contravariant model structures for a left fibration $f$. For the marked case we still need to treat the remaining cartesian anodynes; the hardest case turns out to be the inner anodynes. This case will follow from the following result:

\begin{theorem}
\label{thm:flat}
For $T \rightarrow S$ a categorical equivalence of simplicial sets and $X \rightarrow S$ a left fibration, the morphism $T \times_{S} X \rightarrow X$ in $\sSet$ is a categorical equivalence of simplicial sets. 
\end{theorem}

The theorem is proved by Lurie in \cite[Proposition 3.3.1.3]{HTT}, but his arguments depend on the straightening-unstraightening equivalence. This dependence was removed by Cisinski, who treated the result as \cite[Proposition 5.3.5]{Cisinskibook}. We will shall give a slightly streamlined account of his argument in Section \ref{sec:basechange} as well, see Theorem \ref{thm:flatrestate} specifically.

\subsection{The straightening construction}

We recall the construction of the straightening functor and a number of its basic properties \cite[Section 3.2.1]{HTT}. Given a map of simplicial sets $p \colon X \rightarrow S$, define a simplicial set $S_p$ as the pushout 
\[\begin{tikzcd}
X \ar[r, "i"] \ar[d, "p"] & X * \Delta^0  \ar[d, "q"] \\
            S \ar[r, "j"] & S_p.
\end{tikzcd}\]
Then Lurie defines the \emph{right straightening of $p$} as the functor $\F_{\CC(S_p)}(-,*) \colon \CC(S)^\op \rightarrow \sSet$, which gives us a functor 
\[\Strun \colon \sSet/S \longrightarrow \Fun^\simp(\CC(S)^\op, \sSet).\]
Before we move on to incorporate markings, we record:

\begin{proposition}\label{prop:strcolim}
The functor $\Strun$ preserves colimits and is compatible with base change, in the sense that for every map $p'\colon Y \rightarrow X$ there is a canonical natural isomorphism
\[\Strun(p \circ p') \cong \CC(p)^\op_!\Strun(p').\]
\end{proposition}

We include a short proof in Section \ref{sec:Necklaces}, specifically in Proposition \ref{prop:strcolimrestate}, as an explicit verification seems to be missing from the literature. The proposition formally implies that $\Strun$ has a right adjoint $\Unun \colon \Fun^\simp(\CC(S)^\op,\sSet) \rightarrow \sSet/S$, which is given by sending $G \colon \CC(S)^\op \rightarrow \sSet$ to
\[(p \colon \Delta^n \rightarrow S) \longmapsto \Nat^s(\Strun(\mathrm{id}_{\Delta^n}),\CC(p)^*G)_0.\]
Now let us review how to incorporate markings. Given a map of marked simplicial sets $p \colon X \rightarrow S^\sharp$ and a marked edge $f \colon \Delta^1 \rightarrow X$, say $f \colon x \rightarrow x'$, consider the functor
\[\CC(\Delta^2) = \CC(\Delta^1 * \Delta^0) \xrightarrow{\CC(f * \mathrm{id})} \CC(X * \Delta^0) \longrightarrow \CC(S_p).\]
Now recall that $\F_{\CC(\Delta^2)}(0,2) \cong \Delta^1$. This edge is taken to an edge of $\F_{\CC(S_p)}(px,*) = (\Strun(p))(px)$, which we mark, along with all its images under 
\[(\Strun(p))(px) \xrightarrow{\sigma^*} (\Strun(p))(s)\]
for $\sigma \in \F_{\CC(S)}(s,px)_1$ and $s \in S$ arbitrary. Let us call the resulting functor
\[\Str \colon \mSet/S^\sharp \longrightarrow \Fun^\simp(\CC(S)^\op,\mSet).\]
It acquires a right adjoint by putting the following marking on $\Unun$: from the definition we find for every $f \colon \Delta^1 \rightarrow S$, and $G \colon \CC(S)^\op \rightarrow \mSet$
\[\Hom_{\sSet/S}(\Delta^1, \Unun(G)) = \Nat^\simp(\Strun(f),G)_0\]
and declare an edge of $\Unun(G)$ marked if it corresponds to a transformation that preserves markings when regarded as a transformation $\Str(f) \Rightarrow G$.

The following statements are now simple to check:

\begin{observation}\label{obs:allaboutun}
For $p \colon X \rightarrow S$ and $p' \colon Y \rightarrow X$ one has
\[\Str(p \circ p') \cong \CC(p)_!^\op\Str(p'),\]
and this implies, by adjunction,
\[\Un(G) \times_{S^\sharp} T^\sharp \cong \Un(G \circ \CC(f)^\op)\]
for every map $f \colon T \rightarrow S$ of simplicial sets and every simplicial functor $G\colon \CC(S)^{\mathrm{op}} \to \mSet$. In particular, the fibre of $\Un(G)$ over some $s \in S_0$ is given by $\Un(G(s))$, where $G(s)$ is regarded as a functor $\CC(\Delta^0) = *\rightarrow \mSet$.
\end{observation}

Let us quote from \cite[Proposition 3.2.1.11]{HTT} that $\Str$ takes marked cofibrations and cartesian anodynes to projective cofibrations and projective trivial cofibrations, respectively; see Section \ref{subsec:proj} for a brief review of the projective model structure on functor categories. As observed at the beginning of Section \ref{sec:basechangecocart} this immediately implies:

\begin{theorem}
For every simplicial set $S$ the pair
\[
\begin{tikzcd}
\Str \colon \mSet/S^\sharp \ar[shift left]{r} & \Fun^\simp(\CC(S)^\op,\mSet) \ar[shift left]{l} \cocolon {\Un}
\end{tikzcd}
\] 
is a Quillen adjunction for the cartesian model structure on the left and the projective model structure based on the marked Joyal model structure on the right.
\end{theorem}

The most basic case of straightening is the `absolute' context $S = \Delta^0$, which is straightforward to work out. Without markings one finds
\[
\begin{tikzcd}
{\Strun = |-|_Q} \colon \sSet \ar[shift left]{r} & \sSet \ar[shift left]{l} \cocolon {\mathrm{Sing}_Q = \Unun}
\end{tikzcd}
\] 
with $Q$ the cosimplicial simplicial set described in Section \ref{subsec:coherentnerve}: the simplicial set $Q^n$ is precisely the mapping space from $0$ to $1$ in the simplicial category $\CC((\Delta^{n}*\Delta^0)/\Delta^n)$. Unwinding definitions it can be described as the following quotient of the nerve of $\{S \subseteq [n] \mid S \neq \emptyset\}$:
Two chains $S_0 \subseteq \dots S_k$ and $S'_0 \subseteq \dots \subseteq S'_k$ are identified if there exists an $i \in S \cap S'$ and $S_j \cap [i,n] = S'_j \cap [i,n]$ for all $i$, see e.g.\ \cite[Lemma 3.10]{wurst}. It comes equipped with a natural transformation $Q \Rightarrow \Delta$ induced by taking a vertex $S \subseteq [n]$ in $Q_n$ to $\mathrm{max}(S) \in [n]$, regarded as a vertex in $\Delta^n$. Both $Q$ and $\Delta$ are Reedy cofibrant cosimplicial objects with weakly contractible terms (\cite{wurst} again has detailed verifications), whence Reedy's lemma implies that the induced natural transformation $|-|_Q \rightarrow \mathrm{id}_{\sSet}$ is a weak homotopy equivalence. We also need the more refined statement including markings:

\begin{proposition}\label{cor:fiberUn}
The natural map $\Str(X) \rightarrow X$ just described is a marked categorical equivalence for every $X \in \mSet$. In particular, it induces a natural marked categorical equivalence
\[G(s) \longrightarrow \Un(G)_{s}\]
for every $S \in \sSet$, $s \in S_0$, and projectively fibrant simplicial functor $G \colon \CC(S)^\op \rightarrow \mSet$ and if $S = \Delta^0$ then the adjoint pair $(\Str,\Un)$ is a Quillen equivalence.
\end{proposition}

\begin{proof}
The first statement is \cite[Proposition 3.2.1.14]{HTT} (and the proof there is does not rely on the straightening equivalence). It follows immediately that the functor $\Str \colon \mSet \rightarrow \mSet$ is a left Quillen equivalence for the marked Joyal model structure, which formally implies that the adjoint transformation $\mathrm{id} \rightarrow \Un$ gives a marked categorical equivalence when evaluated on fibrant objects. The second claim then follows from the computation of fibres of unstraightenings in Observation \ref{obs:allaboutun}.
\end{proof}

\begin{corollary}\label{cor:Strprod}
Let $p \colon X \rightarrow S^\sharp$ be a constant map of marked simplicial sets with value $s \in S_0$, which we regard as a map $\Delta^0 \rightarrow S$. Then
\[\Str(p) \cong \CC(s)^\op_!|X|_Q \cong \F_{\CC(S)}(-,s) \times |X|_Q\]
and the natural map to
\[\Str(s) \times X = \F_{\CC(S \cup_{\Delta^0}\Delta^1)}(-,1) \times X\]
is a weak equivalence in $\Fun^\simp(\CC(S)^\op,\mSet)$.
\end{corollary}

\begin{proof}
It only remains to check that composition with the unique edge $s=0 \rightarrow 1$ in $\CC(S \cup_{\Delta^0}\Delta^1)$ induces a categorical equivalence 
\[\F_{\CC(S)}(-,s) \longrightarrow \F_{\CC(S \cup_{\Delta^0}\Delta^1)}(-,1)\]
but it in fact induces an isomorphism by a direct check of the universal property for the category defined this prescription.

The corollary is also a very special case of \cite[Corollary 3.1.1.15]{HTT}.
\end{proof}

\subsection{Functors preserving homotopy colimits}

Let $I$ be a small category and $\C$ a model category (or any category with sufficient colimits equipped with a suitable notion of weak equivalence). Then we can form the homotopy categories $\mathrm{Ho}(\C)$ and $\mathrm{Ho}(\Fun(I,\C))$ by formally inverting the weak equivalences. The \emph{homotopy colimit functor}, if it exists, is the total left derived functor of the colimit functor, i.e., it is a right Kan extension of the composite $\varepsilon \circ \colim_I$ as indicated in the following diagram:
\[
\begin{tikzcd}
\Fun(I,\C) \ar{r}{\colim_I}\ar{d}{\varepsilon} & \C \ar{d}{\varepsilon} \\
\mathrm{Ho}(\Fun(I,\C)) \ar[dashed]{r}{\hocolim_I} & \mathrm{Ho}(\C).
\end{tikzcd}
\]
If $\C$ is a cofibrantly generated model category, then the category $\Fun(I,\C)$ of $I$-diagrams in $\C$ admits the projective model structure. The colimit functor $\colim_I\colon \Fun(I,\C) \rightarrow \C$ is then left Quillen, so that there is a standard recipe for constructing homotopy colimits: given a diagram $X\colon I \rightarrow \C$, one chooses a projectively cofibrant replacement $Y$ of $X$ and $\colim_I Y$ is a representative for $\hocolim_I X$, see also \cite[Chapter IV]{DHKS} or \cite[Part I]{riehl} for textbook accounts of the general theory. In the following definition we write $I^{\triangleright}$ for the right cone on $I$, obtained from $I$ by adjoining a terminal object (which we label $\infty$). Thus an $I^{\triangleright}$-indexed diagram is the same as a cocone on an $I$-shaped diagram. 

\begin{definition}
Let $I$ be a small category and $\C$ some category equipped with a notion of weak equivalences such that the homotopy colimit functor above exists. Then we say that a diagram $X\colon I^{\triangleright} \rightarrow \C$ is a \emph{homotopy colimit diagram}
if the composite 
\[\hocolim_I X \longrightarrow \colim_I X \longrightarrow X(\infty)\]
is an isomorphism in $\mathrm{Ho}(\C)$. We say that a functor between two such categories $\C$ and $\mathcal D$ \emph{preserves homotopy colimits} (of shape $I$) if it preserves all such homotopy colimit diagrams (of shape $I$).
\end{definition}

We will need the following two simple observations:
\begin{itemize}
\item[(1)] Suppose $F\colon \mathcal{C} \rightarrow \mathcal{D}$ is a left Quillen functor between model categories and $Q$ is a functorial cofibrant replacement on $\mathcal{C}$, then $F \circ Q$ preserves homotopy colimits in the sense above: Indeed, this follows easily from the fact that $F$ preserves projectively cofibrant diagrams and weak equivalences between cofibrant objects (by Brown's lemma).
\item[(2)] If two functors $F$ and $G$ are weakly equivalent, then $F$ preserves homotopy colimits if and only if $G$ does.
\end{itemize}

Generally, if $F$ is a left or right Quillen functor we shall write $\mathbf{L}F$ and $\mathbf{R}F$ for the precomposition of $F$ with a (a choice of) functorial cofibrant or fibrant replacement (which exists in all examples we have to consider), respectively. Such composites induce the total left or right derived functor on homotopy categories, but it will be convenient to consider them in this slightly ill-defined form; all statements we shall make about them are easily checked to be invariant under weak equivalences and thus independent of the choices (since any two functorial cofibrant replacements, say, are connected by a zig-zag $Q \Leftarrow Q\circ Q' \Rightarrow Q'$ of weak equivalences). For example, the discussion above is summed up by saying that $\mathbf LF$ preserves homotopy colimits.

\section{Proof of the (un)straightening equivalence}
\label{sec:proof}

In this section we will prove our main result; we restate it for the reader's convenience:

\begin{theorem}
\label{thm:straightening}
For any $\infty$-category $S$ the marked straightening-unstraightening adjunction
\[
\begin{tikzcd}
\Str \colon \mSet/S^\sharp \ar[shift left]{r} & \mathrm{Fun}^\simp(\mathfrak{C}(S)^{\mathrm{op}},\mSet) \ar[shift left]{l} \cocolon \Un
\end{tikzcd}
\]
is a Quillen equivalence for the cartesian model structure on the left, and the projective model structure based on the marked Joyal model structure on the right. 
\end{theorem}

The case of a general simplicial set $T$ (in place of $S$) is not hard to deduce from the statement, as we shall explain at the end of this section. So fix from here on an $\infty$-category $S$. We begin the proof with the following straightforward observation:

\begin{proposition}\label{prop:uncons}
The unstraightening functor
\begin{equation*}
\Un \colon \Fun^\simp(\CC(S)^\op, \mSet) \longrightarrow \mSet/S^\sharp
\end{equation*}
detects weak equivalences between fibrant objects.
\end{proposition}
\begin{proof}
Suppose $\alpha\colon F \rightarrow G$ is a natural transformation between fibrant simplicial functors $\CC(S)^\op\rightarrow \mSet$ such that $\Un(\alpha)$ is a weak equivalence between the cartesian fibrations $\Un(F)$ and $\Un(G)$. In particular, for any $s \in S_0$ the map of fibres
\begin{equation*}
\{s\} \times_{S^\sharp} \Un(F) \longrightarrow \{s\} \times_{S^\sharp} \Un(G)
\end{equation*}
is a categorical equivalence of marked simplicial sets. But then the same is true of $F(s) \rightarrow G(s)$ by Proposition  \ref{cor:fiberUn}.
\end{proof}
The previous proposition implies that $\mathbf{R}\Un$ detects weak equivalences. To prove that the pair $(\Str,\Un)$ is a Quillen equivalence, it therefore suffices to prove that for any map $p\colon X \rightarrow S^\sharp$, the derived unit $\eta$
\begin{equation*}
p \longrightarrow \Un(\Str(p)) \longrightarrow \mathbf{R}\Un(\mathbf L\Str(p))
\end{equation*}
is a cartesian equivalence: From the triangle identity we find that the composite
\[\mathbf{R}\Un F \xrightarrow{\eta}\mathbf{R}\Un \mathbf L\Str \mathbf{R}\Un F \xrightarrow{\mathbf R\Un\epsilon} \mathbf{R}\Un F\]
is (equivalent to) the identity for every $F \colon \CC(S)^\op \rightarrow \mSet$, where $\epsilon$ is the derived counit. But $\eta$ being an equivalence then implies that $\mathbf R\Un\epsilon$ is an equivalence, so \ref{prop:uncons} implies that also $\epsilon$ is one. Note also that $\mathbf L\Str \simeq \Str$ since every object of $\mSet/S^\sharp$ is cofibrant. 

We first prove that the derived unit is an equivalence in the basic case where $X$ is a vertex:

\begin{theorem}
\label{thm:straightenvertex}
For any vertex $s\colon \Delta^0 \rightarrow S$, the derived unit
\begin{equation*}
s \longrightarrow \mathbf{R}\Un(\Str(s))
\end{equation*}
is a cartesian equivalence.
\end{theorem}
\begin{proof}
First observe that
\begin{equation*}
\Str(s)  = \F_{\CC(S)}(-,s)^\sharp.
\end{equation*}
Write $E\F_{\CC(S)}(-,s)^\sharp$ for a fibrant replacement of this functor in $\Fun(\CC(S)^\op, \mSet)$. Then the derived unit may be taken to be the top arrow in the following diagram, picking out the identity morphism of $s$ (or rather its image in the chosen fibrant replacement):
\[
\begin{tikzcd}
s \ar{d}[swap]{\mathrm{id}_s}\ar{r}{\eta} & \Un(E\F_{\CC(S)}(-,s)^\sharp). \\
S^\sharp_{/s} \ar{ur}[swap]{\widetilde{\eta}} &
\end{tikzcd}
\]
Here the vertical map is marked right anodyne (hence a cartesian equivalence) and $\widetilde{\eta}$ is defined as follows. It is the composite of two maps
\begin{equation*}
S^\sharp_{/s} \xrightarrow{u} \Un(\F_{\CC(S)}(-,s)^\sharp) \longrightarrow \Un(E\F_{\CC(S)}(-,s)^\sharp)
\end{equation*}
with the second one being the chosen fibrant replacement. To describe $u$ consider an $n$-simplex $\xi\colon \Delta^n \rightarrow S^\sharp_{/s}$ given by a map
\begin{equation*}
\sigma\colon (\Delta^{n})^\triangleright = \Delta^{n+1} \longrightarrow S
\end{equation*}
sending the final vertex $v$ of $(\Delta^{n})^\triangleright$ to $s$. Defining the relevant map
\begin{equation*}
\xi \longrightarrow \Un(\F_{\CC(S)}(-,s)^\sharp)
\end{equation*}
over $S$ is equivalent to defining a map
\begin{equation*}
\Delta^n \longrightarrow \xi^*\Un(\F_{\CC(S)}(-,s)^\sharp)
\end{equation*}
over $\Delta^n$. Using Observation \ref{obs:allaboutun} once more, this is the same as supplying the adjunct enriched natural transformation
\[\Str(\mathrm{id}_{\Delta^n}) \Longrightarrow \F_{\CC(S)}(\CC(\xi)(-),s)^\sharp\]
of functors $\CC(\Delta^n)^\op \rightarrow \mSet$. Since $\Str(\mathrm{id}_{\Delta^n}) = \F_{\CC((\Delta^{n})^\triangleright)}(-,v)^\sharp$, this is the same as a natural transformation
\begin{equation*}
\F_{\CC((\Delta^{n})^\triangleright)}(-,v)^\sharp \longrightarrow \F_{\CC(S)}(\CC(\xi)(-),s)^\sharp,
\end{equation*}
which is supplied by the functor
\begin{equation*}
\CC(\sigma)\colon \CC((\Delta^{n})^\triangleright) \longrightarrow \CC(S).
\end{equation*}
This defines $u(\xi)$ and thus $\widetilde \eta$.

It now suffices to check that $\widetilde{\eta}$ is a cartesian equivalence. It is a map between right fibrations, so by \ref{cocartmodelstr} (2) this can be verified by checking that the map of fibres over every vertex $t$ of $S$ is a weak homotopy equivalence (since the fibres are fully marked). This map is precisely
\begin{equation*}
\Hom_S^\mathrm R(t,s) \longrightarrow \mathrm{Sing}_Q(E\F_{\CC(S)}(t,s)),
\end{equation*}
where we have used Observation \ref{obs:allaboutun} to identify the fibre of $\Un$. Now, since $\mathrm{Sing}_Q$ is part of a Quillen equivalence (cf. Proposition \ref{cor:fiberUn}) it suffices to verify that the adjoint map
\begin{equation*}
|\Hom_S^\mathrm R(t,s)|_Q \longrightarrow \F_{\CC(S)}(t,s) \longrightarrow E\F_{\CC(S)}(t,s)
\end{equation*}
is a weak homotopy equivalence of simplicial sets. The second map is a weak homotopy equivalence by definition and the first one is a weak homotopy equivalence by Theorem \ref{thm:duggerspivak}.
\end{proof}

To get from Theorem \ref{thm:straightenvertex} to the general straightening theorem we will set up an induction using homotopy colimits. The crucial observation that makes this work is:

\begin{theorem}
\label{thm:Unhocolim}
The functor
\begin{equation*}
\mathbf{R}\Un \colon \Fun^\simp(\CC(S)^{\mathrm{op}},\mSet) \longrightarrow \mSet/S^\sharp
\end{equation*}
preserves homotopy colimits.
\end{theorem}

Its proof requires the following:

\begin{lemma}
\label{lem:fibershocolim}
Let $S$ be a simplicial set and $s \in S_0$ a vertex. Write 
\begin{equation*}
s^*\colon \mSet/S^\sharp \longrightarrow \mSet
\end{equation*}
for the functor assigning to a map $X \rightarrow S^\sharp$ its fibre $\{s\} \times_{S^\sharp} X$. Then its right derived functor
\begin{equation*}
\mathbf{R}s^*\colon \mSet/S^\sharp \longrightarrow \mSet
\end{equation*}
preserves homotopy colimits.
\end{lemma}
\begin{proof}
Factor $s\colon \Delta^0 \rightarrow S$ as a left anodyne $\Delta^0 \rightarrow R$ followed by a left fibration $p\colon R \rightarrow S$. Let $X^\natural \rightarrow S^\sharp$ be a fibrant object of $\mSet/S^\sharp$, i.e.\ a marked cartesian fibration. Then consider the following diagram of pullback squares:
\[
\begin{tikzcd}
\{s\} \times_{S^\sharp} X^\natural \ar{r}{i}\ar{d} & R^\sharp \times_{S^\sharp} X^\natural \ar{r}\ar{d} & X^\natural \ar{d} \\
\Delta^0 \ar{r} & R^\sharp \ar{r}{p} & S^\sharp
\end{tikzcd}
\]
All of the vertical maps are marked cartesian fibrations and the map $i$ is an equivalence of marked simplicial sets by Theorem \ref{thm:Cartfibmarkedanod}. It follows that $\mathbf{R}s^*$ is equivalent to the composition of functors
\begin{equation*}
\mSet/S^\sharp \xrightarrow{\mathbf{R}p^*} \mSet/R^\sharp \longrightarrow \mSet,
\end{equation*}
with the second arrow the forgetful functor. The functor $p^*$ is left Quillen by Theorem \ref{thm:leftfibrightQuillen} and therefore preserves arbitrary weak equivalences (since every object of $\mSet/R^\sharp$ is cofibrant). Hence $\mathbf{R}p^*$ is weakly equivalent to $p^* \simeq \mathbf Lp^*$, which preserves homotopy colimits. The second functor in the display is also left Quillen by Proposition \ref{prop:basechange} and hence preserves homotopy colimits.
\end{proof}

\begin{remark}
\label{rmk:pullbackleftQ}
The lemma implies that for \emph{any} map $f\colon T \rightarrow S$, the pullback functor $\mathbf{R}f^*$ preserves homotopy colimits. Indeed, since weak equivalences between cartesian fibrations are detected fibrewise it suffices to check this after applying $\mathbf{R}t^*$ for every vertex $t$ of $T$. But then $\mathbf{R}t^*\mathbf{R}f^* = \mathbf{R}(ft)^*$ and one reduces to the statement of the lemma above.
\end{remark}

\begin{proof}[Proof of Theorem \ref{thm:Unhocolim}]
Let $I$ be a small category. Since cartesian equivalences between cartesian fibrations over $S$ are detected fibrewise (cf. Theorem \ref{cocartmodelstr}), it suffices to check that for every vertex $s \in S_0$ the natural transformation
\begin{equation*}
\mathbf{R}s^*\hocolim_I \mathbf{R}\Un \longrightarrow \mathbf{R}s^*\mathbf{R}\Un\hocolim_I
\end{equation*}
is a categorical equivalence. By Lemma \ref{lem:fibershocolim} we may rewrite the left-hand side as $\hocolim_I \mathbf{R}s^* \mathbf{R}\Un$. Now recall that Proposition \ref{cor:fiberUn} provides a weak equivalence
\begin{equation*}
\mathbf{R}s^* \mathbf{R}\Un\simeq \mathbf{R}\mathrm{ev}_s \simeq \mathrm{ev}_s;
\end{equation*}
here $\mathrm{ev}_s$ is the functor evaluating a diagram $F\colon \CC(S)^\op \rightarrow \mSet$ at $s$ and the second equivalence is the observation that it preserves arbitrary weak equivalences. We can now conclude the proof by noting that $\mathrm{ev}_s$ also preserves homotopy colimits; indeed, we need only check for example pushouts and coproducts and in both cases 
\[\mathrm{ev}_s \colon \Fun(I,\Fun^\simp(\CC(S)^\op,\mSet)) \longrightarrow \Fun(I,\mSet)\]
preserves projectively cofibrant objects by direct inspection (since projective cofibrations are in particular pointwise cofibrations), whence the statement follows from its obvious counterpart for honest colimits. Alternatively the statement also follows from the fact that $\mathrm{ev}_s$ is a left Quillen functor with respect to the injective model structure on $\Fun(\CC(S)^\op, \mSet)$, which is Quillen equivalent to the projective one.
\end{proof}

We are now in a position to prove the straightening-unstraightening theorem in general.

\begin{proof}[Proof of Theorem \ref{thm:straightening}]
As we have observed, it suffices to prove that for an arbitrary map $f\colon X \rightarrow S^\sharp$ the derived unit $f \rightarrow \mathbf{R}\Un(\Str(f))$  is a weak equivalence of marked simplicial sets over $S$. Starting from Theorem \ref{thm:straightenvertex}, treating the case of a vertex $s\colon \Delta^0 \rightarrow S^\sharp$, we will work our way up to a general map $p$ of marked simplicial sets by a sequence of inductive steps. To this end observe immediately that the class of maps $f\colon X \rightarrow S^\sharp$ for which the derived unit is an equivalence is closed under homotopy colimits by virtue of Theorem \ref{thm:Unhocolim}.

\emph{Step 1)}. Consider a marked 1-simplex $f\colon (\Delta^1)^\sharp \rightarrow S^\sharp$. Then the inclusion of the last vertex $1\colon \Delta^0 \rightarrow (\Delta^1)^\sharp$ is marked right anodyne. Consider the diagram
\[
\begin{tikzcd}
\Delta^0 \ar{d}{\rotatebox{-90}{\(\sim\)}}[swap]{1}\ar{r} & \mathbf{R}\Un(\Str(f \circ 1)) \ar{d}{\rotatebox{-90}{\(\sim\)}} \\
(\Delta^1)^\sharp \ar{r}\ar{dr}[swap]{f} & \mathbf{R}\Un(\Str(f))\ar{d} \\
& S^\sharp.
\end{tikzcd}
\]
The vertical arrows are cartesian equivalences and the top horizontal arrow is a cartesian equivalence by Theorem \ref{thm:straightenvertex}. Hence the bottom horizontal arrow is a cartesian equivalence as well.

\emph{Step 2)}. Consider an unmarked 1-simplex $f\colon (\Delta^1)^\flat \rightarrow S^\sharp$ which is `vertical', in the sense that $f$ factors as
\begin{equation*}
(\Delta^1)^\flat \longrightarrow \Delta^0 \xrightarrow{s} S^\sharp
\end{equation*}
for some vertex $s$ of $S$. Then by Corollary \ref{cor:Strprod} the natural map $\Str(f) \cong \Str(s) \times (\Delta^1)^\flat$ is an equivalence and the derived unit can be identified with the map
\[
\begin{tikzcd}
(\Delta^1)^\flat \ar{rr}\ar{dr}[swap]{f} & & \mathbf{R}\Un(\Str(s)) \times \mathrm{Sing}_Q(\Delta^1)^\flat \ar{dl} \\
& S^\sharp &
\end{tikzcd}
\]
which is the product of the derived unit of $s\colon \Delta^0 \rightarrow S^\sharp$ and the map induced by the transformation $Q \Rightarrow \Delta$ described before Proposition \ref{cor:fiberUn}. The first is a cartesian equivalence by Theorem \ref{thm:straightenvertex}, the second by Proposition \ref{cor:fiberUn}.

\emph{Step 3)}. Consider a general unmarked 1-simplex $f\colon (\Delta^1)^\flat \rightarrow S^\sharp$. Write $(\Delta^2)^\diamond$ for the 2-simplex with the edge $\{1,2\}$ marked and similarly for its horn $(\Lambda^2_2)^{\diamond}$. Consider the map $\varphi\colon (\Delta^2)^\diamond \rightarrow S^\sharp$ that is the degenerate 2-simplex $s_0 f$; it sends the edge $\{0,1\}$ to a degenerate 1-simplex at the vertex $f(0) \in S_0$. Now consider the following diagram of marked simplicial sets, all regarded as sitting over $S^\sharp$ via the map $\varphi$:
\[
\begin{tikzcd}
(\Delta^1)^\flat \ar{r}{\{0,2\}} & (\Lambda^2_2)^\diamond \ar{d} \\
(\Lambda^2_1)^\diamond \ar{r} & (\Delta^2)^\diamond.
\end{tikzcd}
\]
All three arrows are marked right anodyne: The right hand one by definition, and the horizontal ones as the are pushouts of $1\colon \Delta^0 \rightarrow (\Delta^1)^{\sharp}$ and $(\Lambda^2_1)^\flat \rightarrow (\Delta^2)^\flat$, respectively. Thus, to prove that the derived unit is an equivalence on $f$, it suffices to prove it for $\varphi \colon (\Lambda^2_1)^\diamond \rightarrow S^\sharp$. Now observe that the pushout
\begin{equation*}
(\Lambda^2_1)^\diamond = (\Delta^{\{0,1\}})^\flat \amalg_{\{1\}} (\Delta^{\{1,2\}})^\sharp
\end{equation*}
is also a homotopy pushout (both of the maps involved are cofibrations), so that we may consider $(\Delta^{\{0,1\}})^\flat$, $(\Delta^{\{1,2\}})^\sharp$, and $\{1\}$ separately as observed in the preamble of the proof. But the first case follows by step 2) (since it sits over $S$ vertically), the second by step 1), and the third is Theorem \ref{thm:straightenvertex} again.

\emph{Step 4)}. For an unmarked $n$-simplex $f\colon (\Delta^n)^\flat \rightarrow S^\sharp$, we use that the map
\begin{equation*}
(\Delta^{\{0,1\}})^\flat \cup_{\{1\}} (\Delta^{\{1,2\}})^\flat \cup_{\{2\}} \cdots \cup_{\{n-1\}} (\Delta^{\{n-1,n\}})^\flat \longrightarrow (\Delta^n)^\flat
\end{equation*}
is inner anodyne (and thus also marked right anodyne). The domain is an iterated homotopy pushout of unmarked $1$-simplices and vertices, so it is covered by item (3) above. 

\emph{Step 5)}. A general unmarked simplicial set $X^\flat$ over $S^\sharp$ can be written as a homotopy colimit of unmarked $n$-simplices, for example from its skeletal filtration. In detail, suppose we have proved by induction that the derived unit is an equivalence for simplicial sets $X$ of dimension at most $n-1$. Then we write the $n$-skeleton of $X^\flat$ as a pushout
\[
\begin{tikzcd}
\coprod_{\mathrm{nd}(X_n)} (\partial\Delta^n)^\flat \ar{d}\ar{r} & \mathrm{sk}_{n-1} X^\flat \ar{d} \\
\coprod_{\mathrm{nd}(X_n)} (\Delta^n)^\flat \ar{r} & \mathrm{sk}_n X^\flat.
\end{tikzcd}
\]
The coproduct is over nondegenerate $n$-simplices of $X$. The square is also a homotopy pushout since the left vertical map is a monomorphism. The derived unit is an equivalence on all but the lower right corner by what we have already proved. Therefore it is also an equivalence for the $n$-skeleton of $X$. Having proved this for all $n$, we conclude that the derived unit is also an equivalence on the directed colimit of the maps $\mathrm{sk}_n X^\flat \rightarrow S^\sharp$, which is $X^\flat \rightarrow S^\sharp$ and also a homotopy colimit as its structure maps are again monomorphisms. 

\emph{Step 6)}. Finally, a general marked simplicial set over $S^\sharp$ can be obtained from one of the form $X^\flat \rightarrow S^\sharp$ by forming pushouts along maps of the form $(\Delta^1)^\flat \rightarrow (\Delta^1)^\sharp$, adding the markings where necessary. Since these are cofibrations yet again, and $1$-simplices were already treated in steps 1) and 3) this completes the proof.
\end{proof}

\begin{remark}
\begin{enumerate}
\item To obtain the $\infty$-categorical equivalence 
\[\mathrm{Cart}(\mathcal C) \simeq \Fun(\mathcal C^\op,\mathrm{Cat}_\infty)\]
from Theorem \ref{thm:straightening} one proceeds by taking coherent nerves of the bifibrant objects on both sides. Since all objects are cofibrant in the cartesian model structure, and the unstraightening functor is simplicially enriched (by direct inspection), it induces a map between these coherent nerves, which is an equivalence by \cite[Corollary A.3.1.12]{HTT} and Theorem \ref{thm:straightening}. It then remains to identify the coherent nerves of the categories in \ref{thm:straightening} with those in the $\infty$-categorical statement. Using $\cnerv(\mSet) = \mathrm{Cat}_\infty$, there are tautological functors
\[\cnerv(\mSet/S^\sharp) \longrightarrow \mathrm{Cart}(S) \quad \text{and} \quad \cnerv(\Fun^\simp(\CC(S)^\op,\mSet)) \longrightarrow \Fun(S^\op,\mathrm{Cat}_\infty).\]
The former is an equivalence by direct inspection using the fact that the coherent nerve preserves the homotopy type of mapping complexes (see the discussion before Proposition \ref{cor:fiberUn}), and for the latter this is a special case of \cite[Proposition 4.2.4.4]{HTT}.
\item To obtain the statement of Theorem \ref{thm:straightening} also for a general simplicial set $T$ choose a categorical equivalence $f\colon T \rightarrow S$ with $S$ an $\infty$-category and consider the following diagram:
\[
\begin{tikzcd}
\mSet/T^\sharp \ar{d}{f_!}\ar{r}{\Str} & \Fun^\simp(\CC(T)^\op,\mSet) \ar{d}{\CC(f)_!^\op} \\
\mSet/S^\sharp \ar{r}{\Str} & \Fun^\simp(\CC(S)^\op,\mSet).
\end{tikzcd}
\]
It commutes up to natural isomorphism by Observation \ref{obs:allaboutun} and, moreover, the vertical functors are Quillen equivalences: For the functor $f_!$ on the left Lurie deduces this as a consequence straightening theorem in \cite[Proposition 3.3.1.1]{HTT}, but we give a direct argument in Section \ref{subsec:homotopyinvariancecart} below. For the right hand side the statement follows from standard properties of $\CC$ and projective model structures, but as this is strewn over the literature we also include a proof in Section \ref{subsec:proj}. 
\end{enumerate}
\end{remark}

\section{Basic properties of the projective and cartesian model structures}
\label{sec:modelstructures}

The goal of this section is to discuss some aspects of the model structures we used above in more detail. Nothing is new here; the purpose is rather to guide the reader to simple proofs to the requisite inputs, and to collect arguments that are scattered in the literature. We also point out and fill several minor gaps along the way.

\subsection{The projective model structure}\label{subsec:proj}

In the present section we briefly recall the properties of projective model structures on categories of simplicially enriched functors to a Kan simplicial model category. 

The basic existence statement that we used is the following:

\begin{theorem}
Let $\DD$ be a simplicially enriched category and $\M$ be a Kan simplicial model category, that is cofibrantly generated by a set of morphisms $\mathcal I$ in $\mathcal M$. Then $\Fun^\simp(\DD,\M)$ admits a (unique) model structure in which the fibrations and weak equivalences are defined pointwise.

Its cofibrations are generated by the set of enriched natural transformations
\[\{\mathrm{id} \otimes f \colon A \otimes \F_\DD(d,-) \Rightarrow B \otimes \F_\DD(d,-) \mid d \in \DD,  f \colon A \rightarrow B \in \mathcal I\}\]
and it is Kan simplicial for the natural simplicial enrichment of $\Fun^\simp(\DD,\M)$ arising from those of $\DD$ and $\M$.
\end{theorem}

Proofs of this result (and many variants) abound. A particularly clear account in the non-enriched context is \cite[Section 11.6]{Hirschhorn}, and the proof works verbatim to give the model structure in the present generality. That it is again Kan follows easily from the fourth criterion in \cite[Proposition 9.3.7]{Hirschhorn}, as all occuring terms/conditions are formed/verified pointwise.

As an immediate consequence of the definition of fibrations and weak equivalences we obtain the first half of:

\begin{corollary}\label{cor:projquillen}
If $f \colon \CC \rightarrow \DD$ is a simplicially enriched functor and $\M$ is as above, then 
\[
\begin{tikzcd}
{f_!} \colon \Fun^\simp(\CC,\M) \ar[shift left]{r} & \Fun^\simp(\DD,\M) \ar[shift left]{l} \cocolon f^*
\end{tikzcd}
\] 
is a simplicial Quillen adjunction. If $f$ is a Dwyer-Kan equivalence, the adjunction is a Quillen equivalence.
\end{corollary}

Recall that the enriched left Kan extension functor can be described by $\big(f_!(G)\big)(d) \in \M$ being the coequaliser of the two evident maps
\[
\begin{tikzcd}
 \sum_{c, c' \in \CC} G(c) \otimes \big({\F_{\DD}(f(c'),d) \times \F_{\CC}(c,c')}\big) \ar[shift left]{r}\ar[shift right]{r} & \sum_{c \in \CC} G(c) \otimes {\F_{\DD}(f(c),d)}
\end{tikzcd}
\]  
\begin{proof} 
The original reference for the second part is \cite{DwyerKan-invariance} and the statement also appears for example as \cite[Proposition A.3.3.8]{HTT}. For the reader's convenience we include a fairly elementary proof obtained by recombining some arguments from \cite[Section A.3.3]{HTT}: Essential surjectivity of $f$ implies immediately that $\mathbb Rf^*$ is conservative, so it suffices to check that the derived unit is an equivalence. Since $f^*$ preserves weak equivalences, the derived unit is just the ordinary unit on cofibrant object with no need for fibrant resolution. To see that it is an equivalence assume for a moment that $f$ induces anodyne maps on mapping complexes. Then the unit is even a pointwise trivial cofibration: One easily checks that in this case the class of (enriched) transformations $\eta \colon F \Rightarrow F'$ for which the map $f^*f_! F +_F F' \Rightarrow F'$ is a pointwise trivial cofibration is weakly saturated and contains the generating cofibrations, whence it also contains $\emptyset \Rightarrow F$ for any cofibrant $F$ as desired. The general case now follows by embedding $\CC$ into a simplicially enriched category $\mathfrak E$ whose projection to $*$ is a Dwyer-Kan equivalence and applying what we just proved to the first map in 
\[\CC \longrightarrow \DD \times \mathfrak E \xrightarrow{\mathrm{pr}} \DD\]
and a section of the right hand map; a category $\mathfrak E$ as required can for example be produced as follows: Consider the set $X$ of enriched functors $[1](\partial \Delta^n) \rightarrow [1](\Delta^n)$, where $[1](A)$ is the simplicially enriched category with objects $\{0,1\}$, morphism complex from $0$ to $1$ given by $A$ and no other non-identity morphism. We can then employ the small object argument to factor the map $\CC \rightarrow *$ into a functor $\CC \rightarrow \mathfrak E$ in the weakly saturated class generated by $X$ followed by one which admits left lifts against $X$. The latter statement is evidently equivalent to inducing trivial fibrations on mapping spaces, so $\mathfrak E \rightarrow *$ is a Dwyer-Kan equivalence. To see that the first map induces injections on mapping complexes, it suffices (by inspection) to see that this is the case for every pushout of $[1](\partial \Delta^n) \rightarrow [1](\Delta^n)$. But the pushout against some enriched functor $[1](\partial \Delta^n) \rightarrow \mathfrak A$ (with $0 \mapsto x$ and $1 \mapsto y$, say) can be described explicitly: The category $\mathfrak B$ whose objects are those of $\mathfrak A$ and whose morphism complexes are given by 
\[\F_\mathfrak B(a,b) = \left(\sum_{i \geq 0} F_i(a,b)\right)/\sim,\] 
where $F_0(a,b) = \F_{\mathfrak A}(a,b)$,
\[F_i(a,b) = \F_{\mathfrak A}(y,b) \times (\Delta^n \times \F_{\mathfrak A}(y,x))^{i-1} \times \Delta^n \times \F_{\mathfrak A}(a,x)\]
and the equivalence relation is generated (in each simplicial degree) by identifying all $(i+2)$-tuples containing a simplex from the boundary of some $\Delta^n$-factor with the $i$-tuple, in which that simplex is composed with its neighbours (composition is given by concatenating and composing in the evident fashion). It is easy to check that this enriched category has the universal property of the pushout, and also that each simplex in a morphism complex has a unique representative in which none of the $\Delta^n$-components lie in $\partial \Delta^n$, so that the map $\mathfrak A \rightarrow \mathfrak B$ is injective on morphisms as required.
\end{proof}

\begin{proposition}\label{cor:CCQuillen}
If $f \colon S \rightarrow T$ is a categorical equivalence of simplicial sets, then $\CC(f)$ is a Dwyer-Kan equivalence and consequently
\[\begin{tikzcd}{\CC(f)_!}\colon \Fun^\simp(\CC(S),\mSet) \ar[shift left]{r} & \Fun^\simp(\CC(T),\mSet) \ar[shift left]{l} \cocolon {\CC(f)^*}\end{tikzcd}\]
is a Quillen equivalence.
\end{proposition}

The first statement is of course part of the verification that the adjunction
\[\begin{tikzcd}\sSet \ar[shift left]{r}{\CC} & \Cat^\sSet \ar[shift left]{l}{\cnerv}\end{tikzcd}\]
is Quillen for the Joyal and Bergner model structures (and thus a Quillen equivalence by \ref{thm:duggerspivak} and the analogue for $\cnerv$ which is a consequence of the discussion leading up to \ref{cor:fiberUn}). The point of the proof below is that one does not need to even invest the existence of the Bergner model structure.

\begin{proof}
Since Dwyer-Kan equivalences are closed under two-out-of-three, it suffices, by Brown's lemma, to treat the case where $f$ is a trivial cofibration in the Joyal model structure. Choosing an inner anodyne $T \rightarrow T'$ with $T'$ an $\infty$-category further reduces the claim to trivial cofibrations whose target is an $\infty$-category by two-out-of-three. But all such are $J$-anodyne, i.e. lie in the saturation of the inner horn inclusions together with the map $* \rightarrow J$, see e.g. \cite[Proposition 2.5.13]{Land}. Now for these generators $\CC(f)$ has the following property: It admits lifts against all functors $G \colon \mathfrak D \rightarrow \mathfrak E$ of simplicial categories that 
\begin{enumerate}
\item $\mathfrak E$ is Kan enriched
\item $G$ induces Kan fibrations on all mapping spaces (in particular $\mathfrak D$ is also Kan enriched)
\item $G$ induces an isofibration on homotopy categories;
\end{enumerate}
for inner horns this follows from the standard fact that $\CC(\Lambda_i^n) \rightarrow \CC(\Delta^n)$ is an isomorphism away from the mapping complex from $0$ to $n$, where it is an anodyne extension. For $\CC(*) \rightarrow \CC(J)$ it follows for example by adjunction, since $\cnerv(G)$ is an isofibration of $\infty$-categories, which are characterised by lifting against $* \rightarrow J$ as a consequence of Joyal's lifting theorem.

Thus by saturation $\CC(f)$ admits lifts against all such $G$ for an arbitrary $J$-anodyne $f$, and all functors $F \colon \mathfrak A \rightarrow \mathfrak B$ with this lifting property are Dwyer-Kan equivalences by a variant of Quillen's path object argument: For some simplicial category $\mathfrak B$ and some simplicial set $T$, write $E\mathfrak B$ and $F(T,\mathfrak B)$ for the simplicial categories with
\[\F_{E\mathfrak B}(x,y) = E\F_{\mathfrak B}(x,y) \quad \text{and} \quad \F_{F(T,\mathfrak B)}(x,y) = \F_\sSet(T,\F_{\mathfrak B}(x,y)),\]
where $E \colon \sSet \rightarrow \sSet$ is any fibrant replacement in the Kan-Quillen model structure that preserves products (e.g.\ $\mathrm{Ex}^\infty$ or $\mathrm{Sing}|\cdot|$). Then consider the solid part of
\[\begin{tikzcd} \mathfrak A \ar[r,"i"] \ar[d,"F"] & F(\Delta^1,E\mathfrak B) \times_{E\mathfrak B} E\mathfrak A \ar[r,"\mathrm{pr}"] \ar[d,"q"] & E\mathfrak A \\
\mathfrak B \ar[r,"j"] \ar[ru,dotted] & E\mathfrak B\end{tikzcd}\]
where the pullback is formed using $d_0 \colon \Delta^0 \rightarrow \Delta^1$ and $q$ using $d_1$ instead. Then $q$ can also be written as
\[F(\Delta^1,E\mathfrak B) \times_{E\mathfrak B} E\mathfrak A \longrightarrow F(\partial \Delta^1,E\mathfrak B) \times_{E\mathfrak B} E\mathfrak A = E\mathfrak A \times E\mathfrak B \xrightarrow{\mathrm{pr}} E\mathfrak B\]
and therefore satisfies the three conditions above. Thus the dotted lift exists per assumption. But $i$ and $j$ are Dwyer-Kan equivalences (for $i$ this follows from the same statement for the horizontal projection map), and therefore the entire diagram consists of such.
\end{proof}

\subsection{The cartesian model structure}
\label{subsec:cocartesianmodelstruct}

Theorem \ref{cocartmodelstr} describes the cartesian model structure on the category $\mSet/S$, which has the marked injections as its cofibrations and the marked cartesian fibrations $X \rightarrow S$ as its fibrant objects. This model structure was introduced by Lurie in \cite[Section 3.1.3]{HTT} (in the case that $S$ is fully marked) and treated in more generality in \cite[Appendix B]{HA}. Nguyen provides another treatment in \cite{Nguyen} and we have little to add to his elegant arguments; we solely indicate how the dependence on Smith's general existence theorem can be replaced by a more elementary argument, particularly in the case where $S$ is an $\infty$-category. We also elaborate on points (1) and (2) of Theorem \ref{cocartmodelstr}, as well as justify Remark \ref{rmk:cartterminology}.

Recall that we defined the \emph{marked right anodynes} as the weakly saturated class generated by the maps $(A^\marked \rightarrow B^\marked) \boxtimes (C \rightarrow D)$, with $A \rightarrow B$ right anodyne (in $\sSet$) and $C \rightarrow D$ a marked injection. To see that these coincide with Nguyen's \emph{cellular marked right anodynes} (as claimed in Remark \ref{rmk:cartterminology}), one makes the following observation:

\begin{lemma}\label{lem:genmarkright}
The marked right anodynes are generated by the maps
\begin{enumerate}
\item\label{item:markedrightgenerator1} $(\{1\} \rightarrow (\Delta^1)^\sharp) \boxtimes ((\partial \Delta^n)^\flat \rightarrow (\Delta^n)^\flat)$ for $n \geq 0$, and
\item\label{item:markedrightgenerator2} $(\{1\} \rightarrow (\Delta^1)^\sharp) \boxtimes ((\Delta^1)^\flat \rightarrow (\Delta^1)^\sharp)$
\end{enumerate}
and closed under pushout-products with arbitrary marked cofibrations.
\end{lemma}

For the proof we need the following result of Joyal's, see e.g.\ \cite[Proposition 2.1.2.6]{HTT}:

\begin{lemma}
\label{lem:Joyalgenleftanodyne}
The right anodyne maps of simplicial sets are also generated as a weakly saturated class by either of the classes
\begin{enumerate}
\item[(i)] $(\{1\} \rightarrow \Delta^1) \boxtimes (A \rightarrow B)$ with $A \rightarrow B$ running through all injections, or
\item[(ii)] $(\{1\} \rightarrow \Delta^1) \boxtimes (\partial \Delta^n \rightarrow \Delta^n)$ for $n \geq 0$.
\end{enumerate}
\end{lemma}

\begin{proof}[Proof of \ref{lem:genmarkright}]
The closure under pushout-product with marked cofibrations is straightforward. Also, the maps of the lemma are marked right anodyne by definition. To see that they generate all marked right anodynes write $j$ for a cylinder inclusion $(\{1\} \rightarrow \Delta^1) \boxtimes (A \rightarrow B)$ as in the first item of \ref{lem:Joyalgenleftanodyne} and let $k\colon C \rightarrow D$ be a marked injection. Then the pushout-product $j \boxtimes k$ can be rewritten in the form
\[ (\{1\} \longrightarrow (\Delta^1)^\sharp) \boxtimes (i^\sharp \boxtimes k).\]
Since $i^\sharp \boxtimes k$ is a marked injection, this map is clearly in the weakly saturated class generated by (\ref{item:markedrightgenerator1}) and (\ref{item:markedrightgenerator2}).
\end{proof}

Recall also that we defined the \emph{cartesian anodynes} to be saturation of the marked right anodynes, the unmarked inner horn inclusions and $J^\flat \rightarrow (J,0 \rightarrow 1)$. The \emph{marked cartesian fibrations} are those maps which have the right lifting property with respect to the cartesian anodynes. They coincide with Lurie's \emph{marked fibrations} by virtue of the following characterisation:

\begin{proposition}\label{lem:charofmarkedcartfibs}
A map $p \colon E \rightarrow S$ in $\mSet$ is a marked cartesian fibration if and only if
\begin{enumerate}
\item\label{item:charofmarkedcartfibs1} its underlying map is an inner fibration,
\item an edge of $E$ is marked if and only if it is $p$-cartesian and its image in $S$ is marked, and
\item for every marked edge $f\colon t \rightarrow s$ in $S$ and $x \in E$ with $p(x) = s$, there exists a $p$-cartesian lift $g\colon y \rightarrow x$ of $f$.
\end{enumerate}
\end{proposition}

\begin{proof}
This is \cite[Proposition 4.37]{Nguyen}. We shall sketch the proof, in order to comment on a small oversight that is also present in Lurie's version of the result \cite[Proposition 3.1.1.6]{HTT}:

Lifting against unmarked inner anodynes corresponds exactly to item (\ref{item:charofmarkedcartfibs1}), essentially by definition.

Suppose then that $p$ is a marked cartesian fibration. From \cite[Proposition 3.1.1.5]{HTT} we learn that for $n>0$ the inclusion $(\Lambda_n^n,\{n-1 \rightarrow n\}) \rightarrow (\Delta^n,\{n-1 \rightarrow n\})$ is marked right anodyne. The statement implies (with $n=1$) that any marked edge in $S$ with given lift $e$ of its target can be lifted to a marked edge ending at $e$, and (with $n \geq 2$) that every marked edge in $E$ is $p$-cartesian. We are thus left to show that every $p$-cartesian lift of a marked edge in $S$ is marked in $E$. But given such a lift $l \colon d \rightarrow e$, the previous considerations provide a marked edge $k \colon d' \rightarrow e$ with $p(l) = p(k)$. Since $k$ is then also $p$-cartesian, there is an edge $f \colon d \rightarrow d'$, that can be completed to a commutative triangle with $k$ and $l$ (i.e.\ a $2$-simplex in $E$), and because $l$ is $p$-cartesian $f$ is an equivalence in the $\infty$-category $p^{-1}(p(d)) = p^{-1}(p(d'))$. Since this implies that $f$ can be extended to a map $J \rightarrow E$, lifting against $J^\flat \rightarrow (J,\{0 \rightarrow 1\})$ implies that $f$ is marked, and lifting against $(\Delta^2,\{0 \rightarrow 1, 1 \rightarrow 2\})) \rightarrow (\Delta^2)^\sharp$ then forces $l$ to be marked as well. 

Conversely, suppose that $p$ satisfies the conditions of the statement. Then by \cite[Proposition 2.4.1.8]{HTT} $p$ admits lifts against $(\{1\} \rightarrow (\Delta^1)^\sharp) \boxtimes ((\partial \Delta^n)^\flat \rightarrow (\Delta^n)^\flat)$. We next want to argue that it also admits lifts against $(\{1\} \rightarrow (\Delta^1)^\sharp) \boxtimes ((\Delta^1)^\flat \rightarrow (\Delta^1)^\sharp)$. If $S$ is an $\infty$-category this follows from \cite[Proposition 2.4.1.7]{HTT}. In \cite[Proposition 4.37]{Nguyen} Nguyen claims that one can clearly reduce to this situation by pullling $p$ back along the given map $(\Delta^1)^2 \rightarrow S$, and Lurie makes a similar claim regarding $(\Delta^2,\{0 \rightarrow 1, 1 \rightarrow 2\}) \rightarrow (\Delta^2)^\sharp$ in the proof of \cite[Proposition 3.1.1.6]{HTT}. But it is not immediately clear (to us) that this pullback again satisfies the conditions of the proposition: there might a priori be more cartesian than marked edges. For example, Campbell pointed out that it is not generally true that degenerate edges are cartesian \cite{Campbell}, whereas this is true if the base is an $\infty$-category, see \cite[2.1]{rezkdeg} and \cite[Corollary 2.4.1.6]{HTT}. The extra input required is that a locally $p$-cartesian edge (every edge that becomes cartesian in the pullback clearly has this property) is $p$-cartesian, provided it shares its target with a $p$-cartesian lift (even if the base is not assumed an $\infty$-category). Lurie recently recorded this as \cite[Corollary 01U6]{Kerodon}, which allows one to make the desired argument. By \ref{lem:genmarkright} it follows that $p$ admits lifts against all marked right anodynes. Similarly, if $S$ is an $\infty$-category then any map $J \rightarrow E$ represents a $p$-cartesian lift by \cite[Proposition 2.4.1.5]{HTT}, and the general case reduces to this by the same argument. 
\end{proof}

Let us now review the construction of the cartesian model structure: One defines its cofibrations to be the marked cofibrations and its weak equivalences to be those maps
\[\begin{tikzcd}
E \ar[rr,"f"]\ar[rd] && E'\ar[ld] \\
& S &
\end{tikzcd}\]
such that for every marked cartesian fibration $X \rightarrow S$ the induced map 
\[\FHom^\sharp_S(E,X) \longrightarrow \FHom^\sharp_S(E',X)\]
is a weak homotopy equivalence; recall that $\FHom^\sharp_S(E,X) \subset \F_S(E,X)$ denotes the simplicial subset spanned by the vertices representing marked maps (over $S$) and by the edges representing marked maps $(\Delta^1)^\sharp \times E$. Since the cartesian anodynes are closed under pushout-products with marked cofibrations, these are indeed Kan complexes whenever $X \rightarrow S$ is a marked cartesian fibration, see \cite[Proposition 3.1.2.3 \& 3.1.3.6]{HTT} and by varying $X$ the condition is also easily checked equivalent to the a priori weaker requirement that the map above is an equivalence on path components for every marked cartesian fibration $X \rightarrow S$. This is the definition used in \cite[Definition 2.14]{Nguyen}; the cylinder functor $I$ in loc.cit.\ is  $(\Delta^1)^\sharp \times - \colon \mSet/S \rightarrow \mSet/S$ in the case at hand. Finally, the fibrations of the model structure are defined, as they must be, as those maps have the right lifting property against those cofibrations which are simultaneously weak equivalences. Closure under retracts for all three classes of maps and the two-out-of-three property for weak equivalences are then obvious, as is the lifting property for trivial cofibrations against fibrations. For the lifting property for cofibrations against fibrations that are also weak equivalences one first checks that the latter are given precisely by the marked trivial fibrations: Marked trivial fibrations are obviously fibrations in the model structure and that they are weak equivalences is part of \cite[Proposition 2.26]{Nguyen}. The converse now follows by a standard retract trick: Given a fibration $f \colon X \rightarrow Y$ that is also a weak equivalence factor it into a marked cofibration $i \colon X \rightarrow Z$ followed a marked trivial fibration $p \colon Z \rightarrow Y$. By two-out-of-three $i$ is also a weak equivalence whence there exists a filler for the lifting problem
\[\begin{tikzcd}X \ar[r,equal] \ar[d,"i"] & X \ar[d,"f"] \\
 Z \ar[r,"p"] \ar[ru,dotted]& Y.\end{tikzcd}\]
But such a filler expresses $f$ as a retract of $p$, which gives the claim. With this statement in hand it is also clear that every map can be factored into a cofibration followed by fibration that is also a weak equivalence (since the marked cofibration are generated by the set of inclusions $(\partial \Delta^n)^\flat \rightarrow (\Delta^n)^\flat$ and $(\Delta^1)^\flat \rightarrow (\Delta^1)^\sharp$). 

The only difficulty is therefore the existence of factorisations into trivial cofibrations followed by fibrations, since no explicit generators of the trivial cofibrations are known. Both Lurie and Nguyen solve this issue by appeal to Smith's existence theorem, but in the present situation one can be more explicit. Before we do so let us record a few statements that can be obtained without producing such factorisations: As preparation for his proof Nguyen shows in \cite[Lemma 2.31]{Nguyen} that trivial cofibrations with target a marked cartesian fibration are in fact cartesian anodyne, and by \cite[Proposition 2.32]{Nguyen} general trivial cofibrations are characterised by lifting against marked cartesian fibrations $f$ over $S$
\[\begin{tikzcd}
X \ar[rr,"f"]\ar[rd] && Y\ar[ld,"p"] \\
& S &
\end{tikzcd}\]
for which $p$ is also a marked cartesian fibration. This in particular implies that such $f$ are fibrations, in total proving the first addendum to Theorem \ref{cocartmodelstr}. It also shows that the trivial cofibrations form a weakly saturated class and that the cartesian model structure will be simplicial with respect to the enrichment given by $\FHom_S^\sharp$: Using closure of the cartesian anodynes under pushout-products with marked cofibrations, and the consequent fact that $\FHom^\sharp_S(-,X)$ carries marked cofibrations to Kan fibrations for every marked cartesian fibration $X \rightarrow S$, the third condition of \cite[Proposition 9.3.7]{Hirschhorn} is obvious. Finally, one immediately deduces the following statement, compare \cite[Proposition 3.9]{Nguyen}:

\begin{lemma}
\label{lem:cancellation}
Cartesian anodyne morphisms have the right cancellation property among monomorphisms. To be precise, if $i\colon A \rightarrow B$ and $j\colon B \rightarrow C$ are marked injections such that $i$ and $ji$ are cartesian anodyne, then so is $j$.
\end{lemma}

\begin{proof}
The maps $i$ and $ji$ are trivial cofibrations in the cartesian model structure on $\sSet/C$, hence so is $j$. But since the identity $C \rightarrow C$ is a fibrant object, the first item of Theorem \ref{cocartmodelstr} implies that $j$ is cartesian anodyne.
\end{proof}

\begin{remark}
\label{rmk:cancellationrightanodyne}
This entire discussion including Lemma \ref{lem:cancellation} has an evident analogue for the contravariant model structure: A trivial cofibration $A \rightarrow X$ in that model structure with target a right fibration $X \rightarrow S$ is right anodyne and right anodynes consequently have the right cancellation property among monomorphisms. 
\end{remark}

We now return to the problem of factoring an arbitrary map in $\mSet/S$ into a trivial cofibration followed by a fibration in (what we will then have verified to be) the cartesian model structure. For the moment we focus on the case where $S_\flat$ is an $\infty$-category, and discuss the necessary modification for the general case after \ref{lem:countable}. 

\begin{lemma}\label{lem:finitegen}
Replacing $J^\flat \rightarrow (J,0\rightarrow 1)$ in the list of generators for the cartesian anodyne maps by $K^\flat \rightarrow K^\sharp$, where
\[K = \Delta^3/\Delta^{\{0,2\}}, \Delta^{\{1,3\}},\]
yields a class of maps $\mathcal M$ containing the cartesian anodynes, such that for an arbitrarily marked $\infty$-category $S$ any map
\[
\begin{tikzcd}
A \ar{rr}{f}\ar{dr} && B \ar{dl} \\
& S &
\end{tikzcd}
\]
in $\mSet/S$ with $f \in \mathcal M$ is a trivial cofibration.
\end{lemma}

The point of this lemma is that the class $\mathcal M$ is generated by a set of maps between finite marked simplicial sets, which we will use in \ref{lem:countable} below. In fact, we will see that $\mathcal M$ coincides with the class of cartesian anodynes in Corollary \ref{cor:finitegen2}; this statement was first implicitly proven by Lurie as part of \cite[Proposition B.1.6]{HA}. However, he again deduces it as a consequence of the straightening theorem, see item ($\mathrm{A}_1$) in the proof of \cite[Proposition B.1.6]{HA}, and it seems tricky to give a direct argument (cf.\ \ref{rem:fingenannoying} below). We spell out the impact this discussion has on the construction of the cartesian model structure in Remark \ref{rem:fingengeneral}.

\begin{proof}
We note first that the proof of \ref{lem:charofmarkedcartfibs} given above works verbatim upon replacing $J^\flat \rightarrow (J,0 \rightarrow 1)$ by $J^\flat \rightarrow J^\sharp$; indeed, this corresponds to Lurie's choice of generators for the cartesian anodynes in \cite{HTT}.
For the first statement note then only that there is a map $K \rightarrow J$ with $0,2 \mapsto 1$ and $1,3 \mapsto 0$, which expresses the map $J^\flat \rightarrow J^\sharp$ as a pushout of $K^\flat \rightarrow K^\sharp$. 

For the second statement we use the fact that, as explained above, it suffices to check lifting against marked cartesian fibrations $f$
\[
\begin{tikzcd}
X \ar{rr}{f}\ar{dr} && Y \ar[dl,"p"] \\
& S &
\end{tikzcd}
\]
where $p$ is also a marked cartesian fibration. But the assumption on $S$ now implies that both $X_\flat$ and $Y_\flat$ are $\infty$-categories and any map from $K$ to a $\infty$-category takes values in its core. Since equivalences in $X$ are $f$-cartesian, it follows that the marked edge of $K$ goes to a marked edge of $X$ (since this is true in $Y$ by assumption).
\end{proof}

\begin{remark}\label{rem:fingenannoying}
We do not know of a direct presentation of $K^\flat \rightarrow K^\sharp$ in terms of the generators we have chosen for the cartesian anodyne maps (or those from \cite{HTT} or \cite{Nguyen}). Such a presentation exists by virtue of \ref{cor:finitegen2}, and a direct construction would in particular prove Lemma \ref{lem:finitegen} (and also the stronger \ref{cor:finitegen2} itself). Even worse, we do not know of an explicit presentation of the map $J^\flat \rightarrow (J,0 \rightarrow 1)$ in terms of Lurie's or Nguyen's generators (which is why we have chosen it for our exposition), and also do not know whether $K^\flat \rightarrow K^\sharp$ can be replaced by $K^\flat \rightarrow (K,1 \rightarrow 2)$ in the generators for the cartesian anodynes.
\end{remark}

From here the argument for the existence of factorisations (in the case that $S_\flat$ is an $\infty$-category) is essentially the same as for the existence of the Joyal model structure on $\sSet$, e.g.\ in the form given in \cite[Section 8.1]{Dendroids}. It requires the following definition:

\begin{definition}\label{def:deforet}
A \emph{right deformation retract} is a pair of maps of marked simplicial sets
\[
\begin{tikzcd}
A \ar[shift left]{r}{i} & B \ar[shift left]{l}{r}
\end{tikzcd}
\]
with $ri = \mathrm{id}_A$, together with a marked homotopy $h\colon (\Delta^1)^\sharp \times B \rightarrow B$ satisfying $h_0 = \mathrm{id}_B$, $h_1 = ir$, and such that the restriction of $h$ to $(\Delta^1)^\sharp \times A$ is the composition
\begin{equation*}
(\Delta^1)^\sharp \times A \xrightarrow{\pi_A} A \xrightarrow{i} B.
\end{equation*}
\end{definition}

There is of course the dual notion of a left deformation retract given by reversing the direction of the homotopy $h$.

\begin{lemma}
\label{lem:leftdefretractanodyne}
For any right deformation retract as above, the map $i$ is marked right anodyne.
\end{lemma}
\begin{proof}
The diagram
\[
\begin{tikzcd}
A \ar{d}{i}\ar{r}{\{0\} \times \mathrm{id}_A} & (\Delta^1)^\sharp \times A \cup_{\{1\} \times A} \{1\} \times B \ar{d}\ar{r}{\pi_A \cup r} & A \ar{d}{i} \\
B \ar{r}{\{0\} \times \mathrm{id}_B} & (\Delta^1)^\sharp \times B \ar{r}{h} & B
\end{tikzcd}
\]
exhibits $i$ as a retract of the pushout-product of the marked right anodyne map $\{1\} \rightarrow (\Delta^1)^\sharp$ with $i$. Such a pushout-product is again marked right anodyne (cf. Lemma \ref{lem:genmarkright}).
\end{proof}

We shall need the following strengthening of item (1) in Theorem \ref{cocartmodelstr}:

\begin{lemma}\label{lem:deforetract2}
Let $S$ be an arbitrary marked simplicial set. Then a trivial cofibration in $\mSet/S$ between two marked cartesian fibrations is a right (and left) deformation retract over $S$. 
\end{lemma}

\begin{proof}
As explained above a trivial cofibration $i \colon X \rightarrow Y$ as in the statement is in fact cartesian anodyne, so by closure of the cartesian anodynes under pushout products with marked cofibrations the two maps
\[\FHom^\sharp_S(Y,X) \longrightarrow \FHom^\sharp_S(X,X) \quad \text{and} \quad \FHom^\sharp_S(Y,Y) \longrightarrow \FHom^\sharp_S(X,Y)\]
are both trivial fibrations, compare \cite[Remark 3.1.3.4]{HTT}. Lifting the identity of $X$ along the former map provides the retraction $r$ and filling the diagram
\[\begin{tikzcd} \partial \Delta^1 \ar[r,"{(\mathrm{id},ir)}"] \ar[d] & \FHom^\sharp_S(Y,Y) \ar[d] \\
\Delta^1 \ar[r,"\mathrm{const}_i"]\ar[ru,dotted] & \FHom_S^\sharp(X,Y)
\end{tikzcd}\]
provides the deformation.
\end{proof}

\begin{proposition}\label{lem:countable}
For $S$ a marked simplicial set such that $S_\flat$ is an $\infty$-category, the class of trivial cofibrations in $\mSet/S$ is generated by those whose source and target are simplexwise countable, i.e. whose every pullback along some $\Delta^n \rightarrow S_\flat$ is countable.
\end{proposition}

Having found a generating set for the trivial cofibrations we can now apply the small object argument to complete the existence of the cartesian model structure on $\mSet/S$ whenever $S_\flat$ is an $\infty$-category (which is of course the case of main interest).

\begin{proof}
Any trivial cofibration between simplexwise countable sets can be embedded through trivial cofibrations into a trivial cofibration between simplexwise countable marked cartesian fibrations by means of the small object argument applied to the class $\mathcal M$ from Lemma \ref{lem:finitegen}. The latter is then a right deformation retract by the previous lemma. Using these two facts the remainder of the proof is essentially identical to that given in \cite[Lemma 8.13 - Corollary 8.15]{Dendroids} for the analogous statement that the trivial cofibration of the Joyal model structure are generated by those between countable simplicial sets (the same argument is also used in \cite[Section 9.4]{Dendroids} to establish the covariant model structures for dendroidal sets). We therefore refer the reader to loc.cit.\ and refrain from repeating it here.
\end{proof}

Regarding the case of a general base $S$, let us first reiterate that we do not require it for the proof of Theorem \ref{thm:straightening}, since our inductive argument is over a fixed base. To nevertheless obtain the model structure over a general $S$ without appeal to Smith's theorem one can observe, as is also done in \cite[Remark 9.31]{Dendroids}, that Proposition \ref{lem:countable} works mutatis mutandis upon replacing simplexwise countability with a simplexwise bound by any regular cardinal $\alpha$ with the property that a map between objects of $\mSet/S$ that are simplexwise of size smaller than $\alpha$ can be factored into a trivial cofibration followed by a marked cartesian fibration with middle object still simplexwise of size smaller than $\alpha$. If $S_\flat$ is assumed an $\infty$-category Lemma \ref{lem:finitegen} allowed us to choose $\alpha = \aleph_1$ (i.e.\ imposing countability), and since the generating set for the cartesian anodynes we have been using consists of countable objects one can choose $\alpha = \aleph_2$ for the case of arbitrary $S$. We leave the details to the reader. 

\begin{remark}\label{rem:fingengeneral}
Alternatively, one may replace \ref{lem:finitegen} by  \ref{cor:finitegen2} (which states that the class $\mathcal M$ agrees with the cartesian anodynes). This removes assumption that $S_\flat$ be an $\infty$-category from \ref{lem:countable} and so shows that the trivial cofibrations of the cartesian model structure are also generated by those between simplexwise countable objects in general.

Our proof of \ref{cor:finitegen2} below also makes no use of the cartesian model structure (or the straightening equivalence), so no circularity arises.
\end{remark}

Finally, let us discuss item (2) of Theorem \ref{cocartmodelstr}, which concerns the case where the base $S$ is fully marked: The characterisation of fibrant objects as cartesian fibrations with cartesian edges marked is part of Lemma \ref{lem:charofmarkedcartfibs}. The description of weak equivalences between fibrant objects can be seen as a consequence of the following more general result (as we explain below):

\begin{lemma}\label{lem:charcartequiv}
For a diagram
\[\begin{tikzcd}
E \ar[rr,"f"]\ar[rd] && E'\ar[ld] \\
           & S &
\end{tikzcd}\]
of marked simplicial sets whose downwards maps are marked cartesian fibrations, the map $f$ is a cartesian equivalence over $S$ if and only if $f$ induces equivalences of $\infty$-categories 
\[E \times_S {\Delta^1} \longrightarrow E' \times_S {\Delta^1}\]
for all maps $e\colon \Delta^1 \rightarrow S$. 
\end{lemma}
\begin{proof}
This is a special case of \cite[Lemma B.2.4]{HA}, for $\mathfrak P$ the categorial pattern $(\mathcal E, S_2, \emptyset)$ on $S$, where $\mathcal E \subseteq S_1$ is the given marking on $S$. Upon replacing the words `$\mathfrak P$-fibered' and `$\mathfrak P$-anodyne' by `marked cartesian fibration' and `cartesian anodyne', respectively, the proof only uses properties of such maps discussed above, so we refrain from spelling it out here.
\end{proof}

Using the mapping space criterion for cartesian edges from \cite[Proposition 2.4.4.3]{HTT} one readily checks that a map $g$ between cartesian fibrations over an $\infty$-category $T$
	\[\begin{tikzcd}
E \ar[rr,"g"]\ar[rd] && E'\ar[ld] \\
& T &
\end{tikzcd}\]
that preserves cartesian edges is an equivalence if and only if it induces an equivalence on each fibre (using that fully faithful and essentially surjective functors are equivalences). This implies that for every \emph{marked} edge $e \colon x \rightarrow y$ in $S$ the map
\[E \times_S {\Delta^1} \longrightarrow E' \times_S {\Delta^1}\]
is an equivalence if and only if both
\[E \times_S \{x\} \longrightarrow  E' \times_S \{x\} \quad \text{and} \quad E \times_S \{y\} \longrightarrow  E' \times_S \{y\}\]
are equivalences.

To conclude this section let us describe the following special case of the cartesian model structure, which in particular covers the marked Joyal model structure of Lemma \ref{lem:absolutemarked}: 

\begin{lemma}
\label{lem:sliceoverqcat}
If $S$ is an $\infty$-category with its equivalences marked, then the cocartesian and cartesian model structures on $\mSet/S$ both agree with the sliced marked Joyal model structure, which is Quillen equivalent to the sliced (unmarked) Joyal model structure on $\sSet/S_\flat$ using the forgetful functor as the right adjoint.
\end{lemma}

\begin{proof}
It follows from \ref{lem:charofmarkedcartfibs} above (and its dual) that a map $E \rightarrow S$ is a marked cartesian (or cocartesian fibration) if it is an isofibration on underlying simplicial sets and $E$ has precisely its equivalences marked. The first claim now follows for example from the fact that a model structure is determined by its cofibrations and fibrant objects \cite[Proposition E.1.10]{Joyalqcatappl}. 

For the Quillen equivalence to the unmarked Joyal model structure it suffices to treat the case $S=\Delta^0$, since Quillen equivalences pass to slice categories. In this case the left adjoint is given by $(-)^\unmarked$ and the result is straightforward to verify. 
\end{proof}

\section{Base change for the contravariant and cartesian model structures}
\label{sec:basechange}

The purpose of this section is to prove the results promised in Section \ref{sec:basechangecocart}. Specifically, we claimed in Theorem \ref{thm:leftfibrightQuillen} that pullback along a left fibration defines a left Quillen functor with respect to the cartesian model structures; we prove this result in Section \ref{subsec:proofThmbasechange}. As a preliminary, we prove Proposition \ref{thm:Cartfibmarkedanod} in Section \ref{subsec:Cartfibmarkedanod}. In Section \ref{subsec:basechangeinneranod} we will prove Theorem \ref{thm:flat}, stating that pullback along a left fibration preserves categorical equivalences. In Section \ref{subsec:homotopyinvariancecart} we deduce the homotopy invariance of the cartesian model structure quoted at the beginning of Section \ref{sec:proof}.

Most of the material of this section is not original; versions of the arguments in the setting of left and right fibrations can be found in Cisinski's book \cite{Cisinskibook}, and Nguyen extended them to the cartesian case in his thesis \cite{Nguyenthesis}. We hope that the present section provides an efficient path to obtaining these theorems, that could for example be used in a course on $\infty$-categories (as was the motivation for developing our proof of the straightening theorem in the first place).

\subsection{The proof of Proposition \ref{thm:Cartfibmarkedanod}}
\label{subsec:Cartfibmarkedanod}

Recall the statement:

\begin{proposition}\label{thm:Cartfibmarkedanodrestate}
For $T \rightarrow S$ right anodyne and $X \rightarrow S$ a left fibration, the morphism $T \times_{S} X \rightarrow X$ in $\sSet$ is right anodyne. Similarly, for $T \rightarrow S$ marked right anodyne and $X \rightarrow  S$ a marked cocartesian fibration the morphism $T \times_{S} X \rightarrow X$ in $\mSet$ is cartesian anodyne.
\end{proposition}

The marked version was proved by Nguyen \cite[Theorem 4.45]{Nguyen}; the weaker version for right anodynes and left fibrations, see Corollary \ref{cor:rightfibleftanodpullback}, was already observed by Joyal in \cite[Section 11]{Joyalqcatappl}.

To get started recall the notion of a right deformation retract from \ref{def:deforet}.

\begin{lemma}
\label{lem:leftdefretractpullback}
Let $A \xrightarrow{i} B$ be part of a right deformation retract and $X \xrightarrow{p} B$ a marked cocartesian fibration. Then the pullback
\begin{equation*}
A \times_{B} X \longrightarrow X
\end{equation*}
is part of a right deformation retract as well.
\end{lemma}
\begin{proof}
Write $h\colon (\Delta^1)^\sharp \times B \rightarrow B$ for a homotopy exhibiting $A$ as a right deformation retract of $B$. Consider the commutative square
\[
\begin{tikzcd}
(\Delta^1)^\sharp \times (A \times_{B} X) \cup \{0\} \times X \ar{d}\ar{r} & X \ar{d}{p} \\
(\Delta^1)^\sharp \times X \ar{r}{h \circ (\mathrm{id}_{\Delta^1} \times p)} & B.
\end{tikzcd}
\]
The left-hand vertical map is a pushout-product of the marked left anodyne $\{0\} \rightarrow (\Delta^1)^\sharp$ with the monomorphism $A \times_{B} X \rightarrow X$ and hence marked left anodyne. Therefore a lift $k\colon (\Delta^1)^\sharp \times X \rightarrow X$ exists. It is easily verified that $k$ makes $A \times_{B} X$ into a right deformation retract of $X$.
\end{proof}

We can now give the proof of Proposition \ref{thm:Cartfibmarkedanodrestate}, compare \cite[Theorem 4.45]{Nguyen}.

\begin{proof}[Proof of Proposition \ref{thm:Cartfibmarkedanodrestate}]
The class of monomorphisms $T \rightarrow S$ for which the theorem holds is weakly saturated, so it suffices to treat a cylinder inclusion as in the statement of Lemma \ref{lem:Joyalgenleftanodyne}. Consider the diagram of simplicial sets
\[
\begin{tikzcd}
\{1\} \times A \ar{d}[swap]{(1)} \ar{r} & \{1\} \times B \ar{d}[swap]{(2)} \ar{dr}{(3)} & \\
(\Delta^1)^\sharp \times A \ar{r} & \{1\} \times B \cup_{\{1\} \times A} (\Delta^1)^\sharp \times A \ar{r}[swap]{(4)} & (\Delta^1)^\sharp \times B.
\end{tikzcd}
\]
Maps (1) and (3) are clearly part of right deformation retracts, so the conclusion of the proposition holds for them by Lemmas \ref{lem:leftdefretractanodyne} and \ref{lem:leftdefretractpullback}. It also holds for the map (2), since it is a pushout of (1). But then it also holds for (4) by the right cancellation property of cartesian anodynes, cf. Lemma \ref{lem:cancellation}.
\end{proof}

We observe the following consequences of Proposition \ref{thm:Cartfibmarkedanodrestate} for later use:

\begin{corollary}
\label{cor:rightfibleftanodpullback}
If $f\colon X \rightarrow S$ is a left fibration of simplicial sets and $T \rightarrow S$ a right anodyne map, then the pullback
\begin{equation*}
T \times_S X \longrightarrow X
\end{equation*}
is also right anodyne.
\end{corollary}

\begin{proof}
This follows by embedding $\sSet/S$ into $\mSet/S^\sharp$ via the functor $(-)^\sharp$. Alternatively, one observes that the proof above applies verbatim without markings, replacing marked right anodynes (resp. marked cocartesian fibrations) with right anodynes (resp. left fibrations) throughout.
\end{proof}

\begin{corollary}
\label{cor:rightfibleftQuillen}
If $f\colon X \rightarrow Y$ is a left fibration, then
\begin{equation*}
f^*\colon \sSet/Y \longrightarrow \sSet/X
\end{equation*}
is a left Quillen functor with respect to the contravariant model structures on these categories, and so in particular preserves contravariant equivalences.
\end{corollary}
\begin{proof}
Since $f^*$ preserves cofibrations it suffices to observe that it preserves right anodynes (as we explained in the introduction to Section \ref{sec:basechangecocart}), which follows from the previous corollary.
\end{proof}

\subsection{Homotopy invariance of the covariant model structure}
\label{subsec:htpyinvariance}

We write $\mathrm{LFib}(S)$ and $\mathrm{RFib}(S)$ for the homotopy categories of the covariant and contravariant model structure over a simplicial set $S$, respectively. The aim of this section is to prove the following. 
A similar argument as given at the start of Subsection \ref{sec:basechangecocart} shows that 
every map of simplicial sets $f\colon A \rightarrow B$ gives rise to a Quillen adjunction
\[
\begin{tikzcd}
{f_!} \colon \sSet/A \ar[shift left]{r} & \sSet/B \ar[shift left]{l} \cocolon {f^*}.
\end{tikzcd}
\]
Now the goal is to show that this is a Quillen equivalence whenever $f$ is a categorical equivalence. This result and many of the intermediate steps below are originally due to Joyal, our account is a streamlined version Cisinksi's in \cite{Cisinskibook}.

We will repeatedly use the following criterion from \cite[Proposition 4.5.1]{Cisinskibook}:

\begin{proposition}
\label{prop:unitLFib}
Let $f\colon A \rightarrow B$ be a map of simplicial sets and 
\[
\begin{tikzcd}
{\mathbf{L}f_!} \colon \mathrm{LFib}(A) \ar[shift left]{r} & \mathrm{LFib}(B) \ar[shift left]{l}\cocolon{\mathbf{R}f^*}
\end{tikzcd}
\]
the resulting adjunction. Assume that for every vertex $a\colon \Delta^0 \rightarrow A$ the unit
\begin{equation*}
a \longrightarrow \mathbf{R}f^*\mathbf{L}f_!(a)
\end{equation*}
is an isomorphism. Then the unit $\mathrm{id} \rightarrow \mathbf{R}f^*\mathbf{L}f_!$ is an isomorphism for any left fibration $E \rightarrow A$.
\end{proposition}
\begin{proof}
Let $F \rightarrow A$ be any map and pick squares
\[
\begin{tikzcd}
E \ar{r}{i}\ar{d} & \widetilde{E} \ar{d}{p} && F \ar{r}{j}\ar{d} & \widetilde{F} \ar{d}{q} \\
A \ar{r}{f} & B && A \ar{r}{f} & B 
\end{tikzcd}
\]
in which $i$ (resp. $j$) is left (resp. right) anodyne and $p$ (resp. $q$) is a left (resp. right) fibration. Now observe the chain of weak homotopy equivalences of simplicial sets
\begin{eqnarray*}
(\mathbf{R}f^*\mathbf{L}f_! E) \times_A F & \simeq & (\widetilde{E} \times_B A) \times_A F \\
& \cong & \widetilde{E} \times_B F \\
& \xrightarrow{\sim} & \widetilde{E} \times_B \widetilde{F} \\
& \xleftarrow{\sim} & E \times_B \widetilde{F} \\
& \simeq & E \times_A (\mathbf{R}f^*\mathbf{L}f_! F)
\end{eqnarray*}
where the functors on the last line refer to those in the adjunction
\[
\begin{tikzcd}
\mathrm{RFib}(A) \ar[shift left]{r}{\mathbf{L}f_!} & \mathrm{RFib}(B) \ar[shift left]{l}{\mathbf{R}f^*}.
\end{tikzcd}
\]
The arrow on the third (resp. fourth) line is right (resp. left) anodyne by Corollary \ref{cor:rightfibleftanodpullback} (resp. its dual). Now let $a\colon \Delta^0 \rightarrow A$ be an arbitrary vertex and suppose $F \rightarrow A$ is a right fibration. By hypothesis the map $a \rightarrow \mathbf{R}f^*\mathbf{L}f_!(a)$ is a cartesian equivalence over $A$, so by Corollary \ref{cor:rightfibleftQuillen} the map
\begin{equation*}
a \times_A F \longrightarrow \mathbf{R}f^*\mathbf{L}f_!(a) \times_A F
\end{equation*}
is in particular a weak homotopy equivalence of simplicial sets. By the chain of equivalences above we see that also
\begin{equation*}
a \times_A F \longrightarrow a \times_A \mathbf{R}f^*\mathbf{L}f_!(F)
\end{equation*}
is a weak homotopy equivalence. In other words, $F \rightarrow \mathbf{R}f^*\mathbf{L}f_!(F)$ is a fibrewise equivalence between right fibrations and hence a contravariant equivalence over $A$. Applying Corollary \ref{cor:rightfibleftQuillen} we see that for any left fibration $E \rightarrow A$ the map
\begin{equation*}
E \times_A F \longrightarrow E \times_A \mathbf{R}f^*\mathbf{L}f_! F
\end{equation*}
is a weak homotopy equivalence. Once again applying the chain of equivalences at the beginning of the proof, the map
\begin{equation*}
E \times_A F \longrightarrow \mathbf{R}f^*\mathbf{L}f_!(E) \times_A F
\end{equation*}
is also a weak homotopy equivalence. Another application of Corollary \ref{cor:rightfibleftQuillen} shows that his map is even a weak equivalence for arbitrary $F \rightarrow A$. In particular, we may take it to be a vertex inclusion $a\colon \Delta^0 \rightarrow X$, showing that $E \rightarrow \mathbf{R}f^*\mathbf{L}f_!(E)$ is a fibrewise equivalence and hence a covariant equivalence.
\end{proof}

We now prove homotopy invariance of $\mathrm{LFib}$ in the most basic case, compare \cite[Proposition 5.2.1]{Cisinskibook} :

\begin{proposition}
\label{prop:LFibinnerhorn}
For an inner horn inclusion $i\colon \Lambda^n_k \rightarrow \Delta^n$ with $0 < k < n$, the functor
\begin{equation*}
\mathbf{L}i_!\colon \mathrm{LFib}(\Lambda^n_k) \longrightarrow \mathrm{LFib}(\Delta^n)
\end{equation*}
is an equivalence of categories.
\end{proposition}
\begin{proof}
Since $i$ is bijective on vertices, the functor $\mathbf{R}i^*$ is conservative and consequently it suffices to check that the derived unit $\mathrm{id} \rightarrow \mathbf{R}i^*\mathbf{L}i_!$ is an equivalence. By virtue of Proposition \ref{prop:unitLFib} it suffices to do this for the case of a vertex $a\colon \Delta^0 \rightarrow \Lambda^n_k$. Note that $\Delta^n_{a/} \rightarrow \Delta^n$ is a fibrant replacement of $a$ in the covariant model structure over $\Delta^n$, so it suffices to argue that
\begin{equation*}
\Delta^0 \xrightarrow{a} \Lambda^n_k \times_{\Delta^n} \Delta^n_{a/}
\end{equation*}
is a covariant equivalence over $\Lambda^n_k$. For $a > 0$ this is the left anodyne map $a \colon \Delta^0 \rightarrow \Delta^n_{a/}$, whereas for $a = 0$ it is the map $\{0\} \rightarrow \Lambda^n_k$, which is left anodyne: This is clear for $n=2$ and then inductively follows by first recognising both $\{0\} \subseteq \Lambda_i^{n-1} \subseteq \Delta^{\{0,\dots,n-1\}}$ as left anodyne and then filling in the remaining faces of $\Lambda_k^n$ in lexicographical order using e.g.\ \cite[Lemma 1.3.21]{Land} (see also \cite[Lemma 4.4.3]{Cisinskibook} for a different argument).
\end{proof}

The previous proposition in particular implies that any left fibration over the horn $\Lambda^n_k$ extends to a left fibration over the simplex $\Delta^n$ up to covariant equivalence. We will use the following sharper variant, giving extensions up to isomorphism, in order to generalise the result to arbitrary inner anodyne maps:

\begin{corollary}
\label{lem:extendingleftfib}
Let $0 < k < n$. Then for any left fibration $E \rightarrow \Lambda^n_k$, there exists a left fibration $F \rightarrow \Delta^n$ with $E \cong \Lambda^n_k \times_{\Delta^n} F$.
\end{corollary}

In \cite[Theorem 5.2.10]{Cisinskibook} this result is derived from Proposition \ref{prop:LFibinnerhorn} by appeal to the theory of minimal fibrations, which can be entirely bypassed using the following simple argument, which we learned from \cite[Theorem 00ZS]{Kerodon}.

\begin{proof}
Write $i$ for the horn inclusion $\Lambda^n_k \rightarrow \Delta^n$. By Proposition \ref{prop:LFibinnerhorn} there exists a left fibration $p\colon G \rightarrow \Delta^n$ and a covariant trivial cofibration $u\colon E \rightarrow i^*G$ over $\Lambda^n_k$. Define $F \subseteq G$ to be the simplicial subset spanned by those simplices $\sigma\colon \Delta^j \rightarrow G$ such that $i^*\sigma$ factors through $E$. Then $E \rightarrow i^*F$ is an isomorphism; we will prove the lemma by showing that the restriction $p \colon F \rightarrow \Delta^n$ is a left fibration. First of all, note that by (the unmarked version of) Lemma \ref{lem:deforetract2} $u$ is part of a deformation retract. To prove that $F \rightarrow \Delta^n$ is a left fibration, consider a lifting problem
\[
\begin{tikzcd}
A \ar{d}{j}\ar{r}{g} & F \ar[r]\ar{d} & G \ar{dl} \\
B \ar[dashed]{ur}\ar{r}{f} & \Delta^n &
\end{tikzcd}
\]
with $j$ left anodyne. There exists a lift $k \colon B \rightarrow G$ since the slanted map is a left fibration. Now $H$ defines a homotopy (over $\Delta^n$)
\begin{equation*}
h_0 \colon \Delta^1 \times (A \cup_{i^*A} i^*B) = (\Delta^1 \times A) \cup_{\Delta^1 \times i^*A} (\Delta^1 \times i^*B) \xrightarrow{(g \circ \mathrm{pr}_2,H \circ (\mathrm{id} \times i^*k)} G
\end{equation*}
from $g \cup i^*k$ to a map with image in $F$. Now consider the diagram
\[
\begin{tikzcd}
\Delta^1 \times (A \cup_{i^*A} i^*B) \cup \{0\} \times B \ar{r}{h_0 \cup k} \ar{d} & G \ar{d} \\
\Delta^1 \times B \ar{r}{\mathrm{const}_f} & \Delta^n
\end{tikzcd}
\]
The left vertical map is left anodyne. Hence a lift $h$ exists and one verifies that $h_1$ is a solution to the original lifting problem.
\end{proof}

\begin{proposition}
\label{prop:extendingleftfib}
Let $f \colon A \rightarrow B$ be an inner anodyne map of simplicial sets. Then for any left fibration $p\colon E \rightarrow A$ there exists a left fibration $q\colon F \rightarrow B$ with $f^*q = p$.
\end{proposition}
\begin{proof}
Starting from Lemma \ref{lem:extendingleftfib}, it suffices to show that the class of monomorphisms $A \rightarrow B$ having this extension property is closed under pushouts, retracts, (transfinite) compositions. The second and third are completely straightforward. Suppose
\[
\begin{tikzcd}
A \ar{r}{f}\ar{d} & B \ar{d} \\
C \ar{r} & D
\end{tikzcd}
\]
is a pushout square in which $f$ has the desired extension property and let $p\colon E \rightarrow C$ be a left fibration. Take the restriction $A \times_C E$ and extend it to a left fibration $F \rightarrow B$. Then 
\begin{equation*}
E \cup_{E|_A} F \longrightarrow D
\end{equation*}
is easily checked to be a left fibration extending $p$.
\end{proof}

\begin{proposition}
\label{prop:pullbackinneranodyne}
Consider a pullback square
\[
\begin{tikzcd}
X \ar{r}{j}\ar{d} & Y \ar{d}{p} \\
A \ar{r}{i} & B
\end{tikzcd}
\]
with $i$ inner anodyne and $p$ a left fibration. Then $j$ is left anodyne.
\end{proposition}

We will sharpen this in Theorem \ref{thm:pullbackinneranodyne} and show that $j$ is a trivial cofibration in the Joyal model structure.

\begin{proof}
Let $\mathcal{A}$ be the class of maps $i\colon A \rightarrow B$ such that for any left fibration $p\colon Y \rightarrow B$, the pullback of $i$ along $p$ is left anodyne. Then $\mathcal{A}$ is easily seen to be weakly saturated. Hence it suffices to show that it contains the inner horn inclusion $i\colon \Lambda^n_k \rightarrow \Delta^n$. The map $j\colon X \rightarrow Y$ can be interpreted as the derived counit (at $Y$) of the adjunction
\[
\begin{tikzcd}
{\mathbf{L}f_!} \colon \mathrm{LFib}(\Lambda^n_k) \ar[shift left]{r} & \mathrm{LFib}(\Delta^n) \ar[shift left]{l} \cocolon {\mathbf{R}f^*}
\end{tikzcd}
\]
and is thus a covariant equivalence over $\Delta^n$ by Proposition \ref{prop:LFibinnerhorn}. It is also a monomorphism, hence a covariant trivial cofibration, and thus a left anodyne by virtue of Remark \ref{rmk:cancellationrightanodyne}.
\end{proof}

We arrive at the main result of this section, compare \cite[Theorem 5.2.14]{Cisinskibook}.

\begin{theorem}
\label{thm:htpyinvarianceLFib}
If $f\colon A \rightarrow B$ is a categorical equivalence of simplicial sets, then
\begin{equation*}
\mathbf{R}f^*\colon \mathrm{LFib}(B) \longrightarrow \mathrm{LFib}(A)
\end{equation*}
is an equivalence of categories.
\end{theorem}
\begin{proof}
Consider a square
\[
\begin{tikzcd}
A \ar{r}{f}\ar{d} & B \ar{d} \\
\widehat{A} \ar{r} & \widehat{B}
\end{tikzcd}
\]
with vertical maps inner anodyne and $\widehat{A}$ and $\widehat{B}$ both $\infty$-categories. The theorem is straightforward to prove for categorical equivalences between $\infty$-categories; one can use that such an equivalence admits a homotopy inverse, or apply Brown's lemma and reduce to the case of trivial fibrations (in which the necessary verifications are entirely straightforward). By two-out-of-three, it suffices to prove the theorem for the vertical maps. Thus, it suffices to treat the case of an inner anodyne.

For the remainder of this proof, let $f\colon A \rightarrow B$ denote an inner anodyne. First observe that the counit $\mathbf{L}f_!\mathbf{R}f^* \rightarrow \mathrm{id}$ is an isomorphism (so $\mathbf{R}f^*$ is fully faithful): Indeed, for a left fibration $Y \rightarrow B$ the top horizontal map in the square
\[
\begin{tikzcd}
A \times_B Y \ar{r}\ar{d} & Y \ar{d} \\
A \ar{r} & B
\end{tikzcd}
\]
is left anodyne (hence an isomorphism in $\mathrm{LFib}(B)$) by Proposition \ref{prop:pullbackinneranodyne}. It remains to check that $\mathbf{R}f^*$ is essentially surjective, which is immediate from Proposition \ref{prop:extendingleftfib}.
\end{proof}

\subsection{Left fibrations and fully faithful maps}\label{subsec:leftfullyfaith}

A map $f\colon X \rightarrow Y$ between $\infty$-categories is \emph{fully faithful} if for every pair of vertices $x$, $y$ of $X$, the map
\begin{equation*}
\mathrm{Hom}_X(x,y) \longrightarrow \mathrm{Hom}_Y(f(x),f(y))
\end{equation*}
is a weak homotopy equivalence. For a general map $f\colon A \rightarrow B$, we will say it is fully faithful if a fibrant replacement of it is fully faithful, i.e., if there exists a square
\[
\begin{tikzcd}
A \ar{r}{f}\ar{d} & B \ar{d} \\
X \ar{r}{g} & Y
\end{tikzcd}
\]
in which the vertical maps are categorical equivalences, $X$ and $Y$ are $\infty$-categories, and the bottom horizontal map is fully faithful. This condition is easily seen to be independent of the choice of such a square. The point of this short section is the following result, which we will need for Theorem \ref{thm:pullbackinneranodyne} below (compare \cite[Proposition 4.5.2]{Cisinskibook}):

\begin{proposition}
\label{thm:fullyfaithful}
A map $f\colon A \rightarrow B$ of simplicial sets is fully faithful if and only if the functor
\begin{equation*}
\mathbf{L}f_!\colon \mathrm{LFib}(A) \longrightarrow \mathrm{LFib}(B)
\end{equation*}
is fully faithful.
\end{proposition}
\begin{proof}
Using Theorem \ref{thm:htpyinvarianceLFib}, we can take a square as above and without loss of generality assume that $A$ and $B$ are $\infty$-categories. By Proposition \ref{prop:unitLFib} the functor $\mathbf{L}f_!$ is fully faithful if and only if for every vertex $a$ of $A$ the derived unit $a \rightarrow \mathbf{R}f^*\mathbf{L}f_!(a)$ is an isomorphism. This can be made explicit as follows. Observe that $A_{a/} \rightarrow A$ is a fibrant replacement of $a$ and similarly $B_{f(a)/} \rightarrow B$ is a fibrant replacement of $\mathbf{L}f_!(a)$. Hence $\mathbf{L}f_!$ is fully faithful if and only if for every $a$ the map
\begin{equation*}
A_{a/} \longrightarrow A \times_B B_{f(a)/}
\end{equation*}
is a covariant equivalence of left fibrations over $A$. This can be checked on fibres over vertices $b$ of $A$; indeed, the above a map is a weak equivalence if and only if for every $b$ the map
\begin{equation*}
\mathrm{Hom}^\mathrm L_A(a,b) \longrightarrow \mathrm{Hom}^\mathrm L_B(f(a),f(b))
\end{equation*}
is a weak homotopy equivalence of simplicial sets, since the fibres of $A_{a/} \rightarrow A$ identify with the left mapping spaces out of $a$.
\end{proof}

\subsection{Base change of inner anodynes}
\label{subsec:basechangeinneranod}

We can now prove the promised sharpening of Proposition \ref{prop:pullbackinneranodyne}, which will be a key ingredient to proving Theorem \ref{thm:leftfibrightQuillen}:

\begin{theorem}
\label{thm:pullbackinneranodyne}
Consider a pullback square
\[
\begin{tikzcd}
X \ar{r}{j}\ar{d} & Y \ar{d}{p} \\
A \ar{r}{i} & B
\end{tikzcd}
\]
with $i$ inner anodyne and $p$ a left fibration. Then $j$ is $J$-anodyne, and thus in particular a categorical equivalence.
\end{theorem}

Recall that a map of simplicial sets is $J$-anodyne if it lies in the weakly saturated class generated by the inner horn inclusions and $0 \rightarrow J$.

\begin{proof}
Write $\mathcal{A}$ for the class of monomorphisms $i\colon A \rightarrow B$ such that for any left fibration $p\colon X \rightarrow B$ the pullback of $i$ along $p$ is $J$-anodyne. Then $\mathcal{A}$ is weakly saturated, so it suffices to prove that it contains the inner horn inclusions $\Lambda^n_k \rightarrow \Delta^n$ for $0 < k < n$. So consider a pullback square
\[
\begin{tikzcd}
X \ar{r}{j}\ar{d} & Y \ar{d}{p} \\
\Lambda^n_k \ar{r} & \Delta^n.
\end{tikzcd}
\]
It is clear that $j$ is a cofibration. We will argue that it is a categorical equivalence and hence a trivial cofibration in the Joyal model structure; since it has fibrant target it then follows that $j$ is $J$-anodyne. Since $j$ is bijective on vertices, any fibrant replacement of it is essentially surjective; indeed, one can replace a simplicial set by an $\infty$-category without changing its set of vertices. Therefore it suffices to show that $j$ is fully faithful or, by Proposition \ref{thm:fullyfaithful}, that
\begin{equation*}
\mathbf{L}j_!\colon \mathrm{LFib}(X) \longrightarrow \mathrm{LFib}(Y)
\end{equation*}
is fully faithful. This is the same as showing that the unit $\mathrm{id} \rightarrow \mathbf{R}j^*\mathbf{L}j_!$ is an isomorphism. 

To do this, consider any left fibration $q\colon E \rightarrow X$ and form a square
\[
\begin{tikzcd}
E \ar{r}\ar{d}{q} & F \ar{d}{r} \\
X \ar{r}{j} & Y
\end{tikzcd}
\]
in which $E \rightarrow F$ is left anodyne and $r$ is a left fibration. We will conclude the proof by showing that $E \rightarrow X \times_Y F$ is a covariant equivalence over $X$. To this end observe that $X \times_Y F \cong \Lambda^n_k \times_{\Delta^n} F$ and that $E \rightarrow \Lambda^n_k \times_{\Delta^n} F$ is a covariant trivial cofibration over $\Lambda^n_k$ by Proposition \ref{prop:LFibinnerhorn}. But then it is also left anodyne by (the dual of) Remark \ref{rmk:cancellationrightanodyne} and therefore a covariant trivial cofibration over $X$ as well.
\end{proof}

We can now provide the promised proof of Theorem \ref{thm:flat}. We recall the statement:

\begin{theorem}
\label{thm:flatrestate}
For $T \rightarrow S$ a categorical equivalence of simplicial sets and $X \rightarrow S$ a left fibration, the morphism $T \times_{S} X \rightarrow X$ in $\sSet$ is a categorical equivalence of simplicial sets. 
\end{theorem}

\begin{proof}
Let $T \rightarrow S$ be a categorical equivalence of simplicial sets. As in the proof of Theorem \ref{thm:htpyinvarianceLFib}, we may replace $T$ and $S$ by $\infty$-categories via inner anodyne maps, so that by Theorem \ref{thm:pullbackinneranodyne} it suffices to treat the case of a categorical equivalence between $\infty$-categories. Then by Brown's lemma we may further reduce to the case of a trivial fibration between $\infty$-categories. But the pullback of a trivial fibration is a trivial fibration, so in particular a categorical equivalence.
\end{proof}

We close this section by recording the following marked variant of Theorem \ref{thm:pullbackinneranodyne}, which we will need in the next section. Recall that in the absolute case $S = \Delta^0$, the cocartesian and cartesian model structure on $\mSet$ coincide; we dubbed the result the marked Joyal model structure. By Lemma \ref{lem:absolutemarked}, applying $(-)^\flat$ to a categorical equivalence of simplicial sets gives a marked categorical equivalence. Conversely, any marked categorical equivalence between fibrant objects of $\mSet$ gives a categorical equivalence of underlying simplicial sets.

\begin{lemma}
\label{lem:prelimpullbackcocartanod}
Consider a pullback square 
\[
\begin{tikzcd}
X \ar{r}{j}\ar{d} & Y \ar{d}{p} \\
A \ar{r}{i} & B
\end{tikzcd}
\]
of marked simplicial sets with $i$ cartesian anodyne and $p$ a marked cocartesian fibration such that the underlying map of simplicial sets $p_\flat$ is a left fibration. Then $j$ is a marked categorical equivalence.
\end{lemma}
\begin{proof}
As in the proof of Theorem \ref{thm:pullbackinneranodyne}, it suffices to treat the case where $i$ is one of the generators for the class of cartesian anodyne morphisms. If $i$ is marked right anodyne, then the result follows from Proposition \ref{thm:Cartfibmarkedanodrestate}. If $i$ is an (unmarked) inner anodyne, then Theorem \ref{thm:pullbackinneranodyne} implies that the map of simplicial sets $j_{\flat}$ is a trivial cofibration in the Joyal model structure. Since $j$ is a pushout of $(j_\flat)^\flat$, it follows that it is a trivial cofibration in $\mSet$. Finally, if $i$ is the morphism $J^\flat \rightarrow (J, \{0 \rightarrow 1\})$, then the map $X \rightarrow Y$ is an isomorphism of underlying simplicial sets, but just adds some marked edges. To be precise, it is a pushout of the coproduct of maps $(\Delta^1)^\flat \rightarrow (\Delta^1)^\sharp$, with the coproduct ranging over the set of lifts of the edge $\{0 \rightarrow 1\}$ in $J$. Since that edge is an equivalence in the $\infty$-category $J$, any cocartesian lift of it is an equivalence in the $\infty$-category $Y_{\flat}$. Hence $X \rightarrow Y$ can also be written as a pushout of the coproduct of maps $J^\flat \rightarrow (J, \{0 \rightarrow 1\})$ and is therefore cartesian anodyne.
\end{proof}

\subsection{Homotopy invariance of the Cartesian model structure}
\label{subsec:homotopyinvariancecart}

For $S$ a simplicial set, write $\mathrm{Cart}(S)$ and $\mathrm{coCart}(S)$ for the homotopy categories of cartesian and cocartesian fibrations over $S$, respectively. The aim of this section is to prove that $\mathrm{Cart}$ and $\mathrm{coCart}$ take categorical equivalences in $S$ to equivalences of categories, in analogy with the homotopy invariance of $\mathrm{LFib}$ from Theorem \ref{thm:htpyinvarianceLFib}. Much of this works in parallel with the arguments of Section \ref{subsec:htpyinvariance}, but now relying on Proposition \ref{thm:Cartfibmarkedanodrestate} and Theorem \ref{thm:pullbackinneranodyne} rather than the less refined Corollary \ref{cor:rightfibleftQuillen}.

\begin{remark}
The categories $\mathrm{Cart}(S)$ and $\mathrm{coCart}(S)$ are of course the homotopy categories of the cartesian and cocartesian model structures on $\mSet/S^\sharp$, but we will not have to invest the existence of these model structures for the results of this section, so one can use \ref{cor:finitegen2} below to establish their existence over a general base (as explained in \ref{rem:fingengeneral}) without circularity. 

What we shall need to know is that objects of $\mSet/S^\sharp$ can be replaced by (co)cartesian fibrations (to define the functors $\mathbf Lf_!$ below), but this can be achieved by using (co)cartesian anodynes (without investing general factorisations into trivial cofibrations followed by fibrations).
\end{remark}

The following parallels Proposition \ref{prop:unitLFib}:

\begin{proposition}
\label{prop:unitCart}
Let $f\colon A \rightarrow B$ be a map of simplicial sets and 
\[
\begin{tikzcd}
{\mathbf L f_!}\colon \mathrm{coCart}(A) \ar[shift left]{r} & \mathrm{coCart}(B) \ar[shift left]{l} \cocolon {\mathbf{R}f^*}
\end{tikzcd}
\]
the resulting adjunction. If $f$ is fully faithful, then the functor $\mathbf{L}f_!$ is fully faithful.

\end{proposition}
\begin{proof}
We will show that the derived unit $\mathrm{id} \rightarrow \mathbf{R}f^*\mathbf{L}f_!$ is an isomorphism. To do this we mimic the proof of Proposition \ref{prop:unitLFib}. Let $E \rightarrow A^\sharp$ be a marked cocartesian fibration and pick a square
\[
\begin{tikzcd}
E \ar{r}{i}\ar{d} & \widetilde{E} \ar{d}{p} \\
A^\sharp \ar{r}{f} & B^\sharp
\end{tikzcd}
\]
in which $i$ is cocartesian anodyne and $p$ is a marked cocartesian fibration. To show that the map $E \rightarrow \widetilde{E} \times_{B^\sharp} A^\sharp$ is a cocartesian equivalence over $A^\sharp$, it will suffice to show (cf. Theorem \ref{cocartmodelstr}(ii)) that for each vertex $a \in A_0$, the induced map of fibres $E \times_{A^\sharp} a \rightarrow \widetilde{E} \times_{B^{\sharp}} f(a)$ is a categorical equivalence of underlying simplicial sets or, equivalently, that this map is a marked categorical equivalence. Pick a further square of simplicial sets
\[
\begin{tikzcd}
\{a\} \ar{r}{j}\ar{d} & \widetilde{F} \ar{d}{q} \\
A \ar{r}{f} & B 
\end{tikzcd}
\]
in which $j$ is right anodyne and $q$ a right fibration. Then in particular $j^\sharp$ is marked right anodyne and $q^\sharp$ a marked cartesian fibration. The dual of Proposition \ref{thm:fullyfaithful} implies that the map $a \rightarrow A \times_B F$ is a contravariant weak equivalence over $A$, and even a right anodyne by Remark \ref{rmk:cancellationrightanodyne}. Hence $a \rightarrow A^\sharp \times_{B^\sharp} F^\sharp$ is marked right anodyne. Then Proposition \ref{thm:Cartfibmarkedanodrestate} implies that the map
\begin{equation*}
E \times_{A^\sharp} a \rightarrow E \times_{A^\sharp} (A^\sharp \times_{B^\sharp} F^\sharp) = E \times_{B^\sharp} F^\sharp
\end{equation*}
is marked right anodyne as well. Now consider the following commutative diagram:
\[
\begin{tikzcd}
E \times_{A^\sharp} a \ar{r}\ar{d} & \widetilde{E} \times_{B^\sharp} a \ar{d} \\
E \times_{B^\sharp} F^\sharp \ar{r} & \widetilde{E} \times_{B^\sharp} F^\sharp.
\end{tikzcd}
\]
We have just argued that the left vertical map is a marked categorical equivalence; another application of Proposition  \ref{thm:Cartfibmarkedanodrestate} gives the same conclusion for the right vertical map. The dual of Lemma \ref{lem:prelimpullbackcocartanod} implies that the bottom horizontal map is a marked categorical equivalence, using that $ \widetilde{E} \times_{B^\sharp} F^\sharp \rightarrow \widetilde{E}$ is a marked cartesian fibration whose underlying map is a right fibration. (Indeed, it is a pullback of the map $F^\sharp \rightarrow B^\sharp$, which has those two properties.) We conclude that the top horizontal map in the square is a marked categorical equivalence as well, as desired.
\end{proof}

We can now deduce the analogue of Theorem \ref{thm:htpyinvarianceLFib} for (co)cartesian fibrations:

\begin{theorem}
\label{thm:htpyinvarianceCart}
If $f\colon A \rightarrow B$ is a categorical equivalence of simplicial sets, then
\begin{equation*}
\mathbf{L}f_!\colon \mathrm{Cart}(A) \longrightarrow \mathrm{Cart}(B)
\end{equation*}
is an equivalence of categories.
\end{theorem}

\begin{proof}
As in the proof of Theorem \ref{thm:htpyinvarianceLFib} it suffices to treat the case where $f$ is inner anodyne. Then $\mathbf{L}f_!$ is fully faithful by Proposition \ref{prop:unitCart}, so the unit $\eta\colon \mathrm{id} \rightarrow \mathbf{R}f^*\mathbf{L}f_!$ is an isomorphism. The map $f$ is a bijection on vertices; by the fibrewise criterion for weak equivalences between cocartesian fibrations, it follows that the right adjoint $\mathbf{R}f^*$ is a conservative functor. This implies that the counit $\varepsilon\colon \mathbf{L}f_!\mathbf{R}f^* \rightarrow \mathrm{id}$ is an isomorphism as well. Indeed, it suffices to check this after composing with $\mathbf{R}f^*$, where it follows from the triangle identity
\[
\begin{tikzcd}
\mathbf{R}f^* \ar{r}{\cong}\ar[equal]{dr} & \mathbf{R}f^*\mathbf{L}f_!\mathbf{R}f^* \ar{d}{\mathbf{R}f^*\varepsilon} \\
& \mathbf{R}f^*. 
\end{tikzcd}
\]
\end{proof}

Theorem \ref{thm:htpyinvarianceCart} implies that for a categorical equivalence map $A \rightarrow B$ and a cartesian fibration $p\colon X \rightarrow A$, there exists a cartesian fibration over $B$ and an equivalence between its pullback to $A$ and the map $p$. The same argument used to prove Lemma \ref{lem:extendingleftfib} (now with markings) provides the following sharpening:

\begin{proposition}
\label{prop:extendingCartfib}
Let $f \colon A \rightarrow B$ be an inner anodyne map of simplicial sets. Then for any cartesian fibration $p\colon E \rightarrow A$ there exists a cartesian fibration $q\colon F \rightarrow B$ with $f^*q = p$.
\end{proposition}

The decisive case of an inner horn inclusion first appeared as part of \cite[Theorem 3.4.7]{Nguyenthesis}, who developed the notion of minimal cartesian fibrations for the proof, compare the comment after Corollary \ref{lem:extendingleftfib}. The statement is also obtained in \cite[Corollary 029F]{Kerodon}, but the proof there uses (a weak form of) the straightening theorem as input. 

Proposition \ref{prop:extendingCartfib} also implies the following strengthening of Lemma \ref{lem:finitegen}:

\begin{corollary}\label{cor:finitegen2}
One can replace $J^\flat \rightarrow (J,0\rightarrow 1)$ in the list of generators for the cartesian anodyne maps by $K^\flat \rightarrow K^\sharp$, where
\[K = \Delta^3/\Delta^{\{0,2\}}, \Delta^{\{1,3\}}.\]
\end{corollary}

This result removes the assumption that the base be an $\infty$-category from Lemma \ref{lem:countable}, which in turn allows one to establish the cartesian model structure of a general base by using trivial cofibrations between simplexwise countable objects, see Remark \ref{rem:fingengeneral}. An interpretation of the simplicial set $K$ is as follows: it is the universal example of a 1-simplex, namely the edge from 1 to 2, together with 2-simplices witnessing left and right inverses to it.

\begin{proof}
To see that $K^\flat \rightarrow K^\sharp$ is cartesian anodyne, we have to show that any lifting problem against a marked cartesian fibration admits a filler. This immediately reduces to the case where this fibration has target $K^\sharp$ and is therefore a cartesian fibration with its cartesian edges marked. Writing it as the pullback of a cartesian fibration over an $\infty$-category by means of Proposition \ref{prop:extendingCartfib}, then makes the claim obvious since all edges of $K$ go to equivalences in any $\infty$-category (and equivalences are always cartesian, hence marked).

For the converse, note only that one can replace the generator $J^\flat \rightarrow (J,0 \rightarrow 1)$ in the definition of cartesian anodynes by $J^\flat \rightarrow J^\sharp$; the argument for Proposition \ref{lem:charofmarkedcartfibs} goes through verbatim (and in fact, the $J^\flat \rightarrow J^\sharp$ is Lurie's choice in \cite{HTT}). But $J^\flat \rightarrow J^\sharp$ is clearly a pushout of $K^\flat \rightarrow K^\sharp$.
\end{proof}

Another application is the following, which is also deduced from the straightening theorem in \cite{HTT}:

\begin{corollary}
\label{cor:cartfibcatfib}
A cartesian fibration $p\colon E \rightarrow A$ is in particular a categorical fibration (i.e., a fibration in the Joyal model structure).
\end{corollary}
\begin{proof}
If the simplicial set $A$ is an $\infty$-category, this is straightforward: a cartesian fibration $p$ is in particular an isofibration, and a map between $\infty$-categories is a categorical fibration if and only if it is an isofibration, see e.g. \cite[Theorem 2.5.14]{Land}. In the general case, pick an inner anodyne map $A \rightarrow B$ such that $B$ is an $\infty$-category. By Proposition \ref{prop:extendingCartfib} there exists a cartesian fibration $q\colon F \rightarrow B$ with $f^*q = p$. Then $q$ is a categorical fibration by the first sentence, hence so is its pullback $p$.
\end{proof}

\subsection{The proof of Theorem \ref{thm:leftfibrightQuillen}}
\label{subsec:proofThmbasechange}

Recall that the statement:

\begin{theorem}
\label{thm:leftfibrightQuillenrestate}
If $f\colon T \rightarrow S$ is a left fibration, then 
\[
\begin{tikzcd}
{f^*} \colon \mSet/S^\sharp \ar[shift left]{r} & \mSet/T^\sharp \ar[shift left]{l} \cocolon {f_*}
\end{tikzcd}
\]
is a Quillen adjunction with respect to the cartesian model structures.
\end{theorem}

This first appeared as \cite[Proposition 3.6.4]{Nguyenthesis}. As remarked before, it suffices to show that $f^*$ sends cartesian anodynes to trivial cofibrations, because then its right adjoint $f_*$ will preserve fibrations between fibrant objects. For the marked right anodynes and the cartesian anodynes of the form $J^\flat \rightarrow (J, \{0 \rightarrow 1\})$ we already checked this in the proof of Lemma \ref{lem:prelimpullbackcocartanod}. It only remains to show that for a morphism $i$ of the form
\[
\begin{tikzcd}
(\Lambda_k^n)^\flat \ar{rr}{i}\ar{dr} && (\Delta^n)^\flat \ar{dl} \\
& S^\sharp, &
\end{tikzcd}
\]
the pullback of $i$ along $f$ is a trivial cofibration in the cartesian model structure over $T^\sharp$. Observe that this pullback is a pushout of the morphism
\[(T \times_S \Lambda_k^n)^\flat \rightarrow (T \times_S \Delta^n)^\flat.\]
The underlying morphism of simplicial sets is a trivial cofibration in the Joyal model structure by Theorem \ref{thm:pullbackinneranodyne}. Hence, the desired conclusion follows from:

\begin{lemma}
Let $S$ be a simplicial set. The functor 
\[\sSet/S \rightarrow \mSet/S^\sharp\colon (X \rightarrow S) \mapsto (X^\flat \rightarrow S^\sharp) \]
is left Quillen from the sliced Joyal model structure to the cartesian model structure.
\end{lemma}
\begin{proof}
It is clear that the functor of the proposition is left adjoint and preserves cofibrations. Also, its right adjoint (the forgetful functor) preserves fibrant objects by Corollary \ref{cor:cartfibcatfib}. Now let $i\colon A \rightarrow B$ be a trivial cofibration in $\sSet/S$. Then $i^\flat$ is a trivial cofibration in $\mSet/S^\sharp$ if and only if for every marked cartesian $p\colon X^\natural \rightarrow S^\sharp$, the map 
\[\FHom_S^\sharp(B^\flat, X^\natural) \rightarrow \FHom_S^\sharp(A^\flat, X^\natural)\]
is a homotopy equivalence. But this map is identified with the restriction of
\[\mathrm{Fun}_S(B, X) \rightarrow \mathrm{Fun}_S(A, X),\]
to cores, and this map is a trivial fibration because $i$ is a trivial cofibration and $X \rightarrow S$ a categorical fibration.
\end{proof}

\section{Necklaces and path categories}\label{sec:Necklaces}

For a simplicial set $S$, the mapping complexes in the simplicial category $\CC(S)$ admit a rather explicit description in terms of \emph{necklaces} in the simplicial set $S$ by work of Dugger and Spivak \cite{duggerspivakrigid}. In the present section we review this description, use it to verify Proposition \ref{prop:strcolim} and, closely following the arguments of Dugger and Spivak, also Theorem \ref{thm:duggerspivak}.

\subsection{Necklaces and the straightening construction}\label{subsec:straighteningadjunction}

We denote the category $\partial \Delta^1 /\sSet$ of doubly pointed simplicial sets by $\sSet_{\ast,\ast}.$ The marked points of a doubly pointed simplicial set $X$ are labelled $0$ and $1$. For two doubly pointed simplicial sets $X$ and $Y$ we define their wedge $X \vee Y$ by putting them `end to end'. To be precise, it is the pushout $X \cup_{\Delta^0} Y$ where $\Delta^0$ includes as $1$ in $X$ and $0$ in $Y$. This wedge is regarded as an object of $\sSet_{\ast,\ast}$ by marking the point $0$ of $X$ and $1$ of $Y$. We regard the standard $n$-simplex as an object of $\sSet_{\ast,\ast}$ by marking its endpoints $(0,n)$.

\begin{definition}
A \emph{necklace} is a doubly pointed simplicial set of the form
\[ \NL = \Delta^{n_0} \vee \dotsb \vee \Delta^{n_k} \]
with either all $n_i$ strictly larger than $0$ or $k = 0$. The constituent simplices $\Delta^{n_i}$ are called the \emph{beads} of the necklace. The vertices along which the simplices have been glued, together with the initial and terminal vertex, are called the \emph{joints}. For a necklace $\NL$, we will write $\alpha(\NL)$ for its initial vertex and $\omega(\NL)$ for its terminal vertex. We define the category $\Nec$ of necklaces to be the full subcategory of $\sSet_{\ast,\ast}$ spanned by the necklaces.
\end{definition}

The joints of a necklace $\NL$ can be equipped with an evident linear ordering in which $\alpha(\NL)$ is initial and $\omega(\NL)$ terminal. From the defining colimit diagram it is not difficult to see that for two joints $i,j$ of a necklace $\NL$ we have 
\[\F_{\CC(\NL)}(i,j) \cong \prod_{k=i+1}^j \F_{\CC(\Delta^{n_k})}(0,n_k) \cong \prod_{k=i+1}^j \nerve(\mathrm{P}_{0,n_k}),\]
with composition translating to the identity map on the right hand side. Mapping complexes in other path categories can then be described in terms of these as follows:

\begin{proposition}\label{prop:Cnecklace}
For any $S \in \sSet$ there is a canonical isomorphism
\[\F_{\CC(S)}(s,t) \cong \colim_{\NL \rightarrow S} \F_{\CC(\NL)}(\alpha(\NL),\omega(\NL)),\]
where the colimit runs over the full subcategory of the slice category $\Nec / S$ consisting of maps $f$ with $f(\alpha(\NL)) = s$ and $f(\omega(\NL)) = t$.
\end{proposition}
Unwinding definitions, this means that 
\begin{enumerate}
\item[i)] every $k$-simplex on the left can be described by a necklace mapping to $S$, say $f \colon \NL \rightarrow S$, such that $f(\alpha(\NL)) = s$ and $f(\omega(\NL)) = t$, together with a $k$-flag $T_0 \subseteq \dots \subseteq T_k$ of vertices of $\NL$ containing all joints, 
\item[ii)] two such necklaces determine the same vertex if and only if they are equivalent under the equivalence relation generated by flag-preserving maps of necklaces over $S$, and
\item[iii)] the simplicial structure corresponds to removing and doubling elements in the flag.
\end{enumerate}

For the proof, one first checks that the right hand side defines a simplicially enriched category $C(S)$, with composition induced by the concatenation of necklaces taking $M,N \in \Nec$ to the pushout
\[\begin{tikzcd}
\Delta^0 \ar[r,"\alpha(N)"] \ar[d,"\omega(M)"] & N \ar[d] \\
M \ar[r] &                   M \vee N.
\end{tikzcd}\]
One checks that the resulting functor $C \colon \sSet \rightarrow \Catsimp$ satisfies the universal property defining $\CC$. Neither step is difficult and both are carried out in detail in \cite[Section 4]{duggerspivakrigid}, so we refrain from elaborating. Another useful fact is that every equivalence class of flagged necklaces as above contains a unique representative such that 
\begin{enumerate}
\item every bead represents a non-degenerate simplex of $S$, and
\item the flag starts with $T_0$ being the set of joints of $\NL$ and is exhaustive, i.e.\ $T_k = \NL_0$.
\end{enumerate}
Let us call such representatives totally non-degenerate. For the proof, which is again not difficult, see \cite[Corollary 4.8]{duggerspivakrigid}.

Before we dive into the homotopical analysis of mapping complexes in the categories $\CC(S)$, let us use the above description to verify Proposition \ref{prop:strcolim}. We repeat the statement:

\begin{proposition}\label{prop:strcolimrestate}
The functor $\Strun$ preserves colimits and is compatible with base change, in the sense that for maps $p \colon X \rightarrow S$ and $p'\colon Y \rightarrow X$ there is a canonical natural isomorphism
\[\Strun(p \circ p') \cong \CC(p)^\op_!\Strun(p').\]
\end{proposition}


Recall that $\Strun(p) = i^*\F_{\CC(S_p)}(-,*)$ where $S_p$ is defined by the pushout
\[\begin{tikzcd} X \ar[r,"j"] \ar[d,"p"]& X \ast \Delta^0 \ar[d,"q"] \\
                 S \ar[r,"i"] & S_p.\end{tikzcd}\]
For the proof of \ref{prop:strcolimrestate} we need the following lemma:

\begin{lemma}\label{lem:pushoutfaith}
Let
\[\begin{tikzcd} \mathfrak A \ar[r,"j"] \ar[d,"p"] & \mathfrak B \ar[d,"q"] \\
                 \mathfrak C \ar[r,"i"] & \mathfrak D\end{tikzcd}\]
be a pushout diagram of simplicial categories such that $j$ is fully faithful. Then also $i$ is fully faithful and the square
\[\begin{tikzcd} \Fun^\simp(\mathfrak A^\op,\sSet) & \ar[l,"j^*"] \Fun^\simp(\mathfrak B^\op,\sSet) \\
\Fun^\simp(\mathfrak C^\op,\sSet) \ar[u,"p^*"] & \Fun^\simp(\mathfrak D^\op,\sSet) \ar[l,"i^*"] \ar[u,"q^*"] \end{tikzcd}\]
is left adjointable in the sense that
\[q^*i_! \cong j_!p^*\]
via the Beck-Chevalley transformation.
\end{lemma}

Here we mean full faithfulness in the strict sence that $i$ induces isomorphisms (and not just weak homotopy equivalences) 
\[\F_{\mathfrak A}(x,y) \longrightarrow \F_{\mathfrak B}(i(x),i(y)).\]

\begin{proof}
By Yoneda's lemma we have a commutative diagram
\[\begin{tikzcd}\mathfrak C \ar[d,"i"] \ar[r] & \Fun^\simp(\mathfrak C,\sSet) \ar[d,"i_!"] \\
\mathfrak D \ar[r] & \Fun^\simp(\mathfrak D,\sSet)\end{tikzcd}\]
with fully faithful horizontal maps, so that we may as well check $i_!$ fully faithful. We claim that, under the isomorphism
\[\Fun^\simp(\mathfrak D,\sSet) \cong \Fun^\simp(\mathfrak C,\sSet) \times_{\Fun^\simp(\mathfrak A,\sSet)} \Fun^\simp(\mathfrak B,\sSet)\]
induced by the pullback functors, the functor $i_!(F)$, for some $F \colon \mathfrak C^\op \rightarrow \sSet$ corresponds to the pair $(F,j_!p^*F)$, which defines an element in the pullback above, precisely because $j$ is fully faithful so that $j^*j_!$ is the identity of $\Fun^\simp(\mathfrak A,\sSet)$. Granting this claim for a moment, we find that $i^*i_!F = F$ which shows that $\iota_!$ is fully faithful, and also $q^*i_!F = j_!p^*F$ by direct inspection.

To obtain the claim we simply observe 
\[\Nat(i_! F,G) = \Nat(F,i^*G) = \Nat(F,i^*G) \times_{\Nat(p^*F,(qj)^*G)} \Nat(j_!p^*F,q^*G)\]
for $G \colon \mathfrak D^\op \rightarrow \sSet$, whence we are done by another application of Yoneda's lemma.
\end{proof}

\begin{proof}[Proof of Proposition \ref{prop:strcolimrestate}]
The commutation of $\Strun$ with colimits is essentially immediate from the necklace description of mapping complexes in path categories: a $k$-simplex in $(\Strun(p)\big)(s) = \F_{\CC(S_p)}(s,\ast)$ is represented by a unique totally non-degenerate necklace from $s$ to $\ast$ in $S_p$. All but the last bead then lie in $S \subseteq S_p$, and the final bead is necessarily of the form 
\[\Delta^{n} = \Delta^{n-1} * \Delta^0 \xrightarrow{f \ast \mathrm{id}} X \ast \Delta^0,\]
for some unique $f$ and $n$, i.e. it is the same as a non-degenerate simplex in $X$. This data along with a flag clearly depends on $X$ in a colimit preserving manner.

We next prove the universal case of the compatibility of straightening with Kan extensions, i.e.\ that 
\[\Strun(p) \cong \CC(p)_!\CC(j)^*\F_{\CC(X \ast \Delta^0)}(-,*).\]
The general unmarked case then follows since Kan extensions compose, and the marked case by direct inspection. To prevent confusion let us first remark that the statement does \emph{not} directly follow from the adjointability statement of the previous lemma (which goes the wrong way).

The necklace description of $\CC(X \ast \Delta^0)$ does, however, immediately imply that $\F_{\CC(X \ast \Delta^0)}(-,*)$ is almost left Kan extended from $\CC(X)^\op$ in the sense that the natural transformation
\[\CC(j)_!\CC(j)^*\F_{\CC(X \ast \Delta^0)}(-,*) \Longrightarrow \F_{\CC(X \ast \Delta^0)}(-,*)\]
is an isomorphism on all objects but $* \in X \ast \Delta^0$. From the pointwise formula one readily computes that the left hand side evaluates to $\emptyset$ at $\ast$, whereas $\F_{\CC(X \ast \Delta^0)}(*,*) = *$ (this is again obvious from the necklace description for example). Let us denote the two sides by $\mathrm{F}_\emptyset$ and $\mathrm{F}_\ast$, respectively. The claim then becomes, that the Beck-Chevalley map
\[\CC(p)_!\CC(j)^*\mathrm{F}_\ast \Longrightarrow \CC(i)^*\CC(q)_!\mathrm{F}_\ast\]
is an isomorphism. Now Lemma \ref{lem:pushoutfaith} implies that also $\CC(i)$ is fully faithful, so one can compute
\begin{align*}
\CC(i)^*\CC(q)_!\mathrm{F}_\emptyset &\cong \CC(i)^*\CC(q)_!\CC(j)_!\CC(j)^*\mathrm{F}_\ast \\
&\cong \CC(i)^*\CC(i)_!\CC(p)_!\CC(j)^*\mathrm{F}_\ast \\
&\cong \CC(p)_!\CC(j)^*\mathrm{F}_\ast
\end{align*}
We will therefore be done if we can show that the natural map
\begin{equation}\tag{$\ast$}\label{eq:middleofaproof}
\big(\CC(q)_!\mathrm{F}_\emptyset\big)(s) \longrightarrow \big(\CC(q)_!\mathrm{F}_\ast\big)(s)
\end{equation}
is an isomorphism for every $s \in S$. Using the pointwise formula for enriched Kan extensions (we recalled it after \ref{cor:projquillen}) the two sides evaluate to the coequaliser of the two composition maps
\[\sum_{c,c' \in X \cup \{\ast\}} \mathrm{F}_\epsilon(c) \times \F_{\CC(S_p)}(s,q(c')) \times \F_{\CC(X \ast \Delta^0)}(c',c) \longrightarrow \sum_{c \in X \cup \{\ast\}} \mathrm{F}_\epsilon(c) \times \F_{\CC(S_p)}(s,q(c))\]
for $\epsilon \in \{\emptyset,\ast\}$. These terms differ only in the summands with $c= \ast$ which are empty for $\epsilon = \emptyset$. But this change has no influence on the coequaliser: To see surjectivity, observe that the extraneous summand on the right evaluates to $\mathrm{F}_\ast(\ast) \times \F_{\CC(S_p)}(s,*)$ with $\mathrm{F}_\ast(\ast)= \ast$. But by the necklace description of $\F_{\CC(S_p)}(s,*)$ any simplex $\sigma$ in it lies in the image of the composition map
\[\mathrm{F}_\ast(\ast) \times \F_{\CC(S_p)}(s,q(c')) \times \F_{\CC(X \ast \Delta^0)}(c',*) \longrightarrow \mathrm{F}_\ast(\ast) \times  \F_{\CC(S_p)}(s,*),\]
where $c' \in X$ is a preimage of the second to last joint in a flagged necklace representing $\sigma$ (such a preimage exists since the entire last bead of a necklace with target $* \in S_p$ lifts to $X \ast \Delta^0$). Applying the other composition map
\[\mathrm{F}_\ast(\ast) \times \F_{\CC(S_p)}(s,q(c')) \times \F_{\CC(X \ast \Delta^0)}(c',*) \longrightarrow \mathrm{F}_{\ast}(c') \times \F_{\CC(S_p)}(s,q(c'))\]
to such a preimage then provides an element which 
\begin{enumerate}
\item by definition represents the same element as $\sigma$ in the coequaliser for $\epsilon = \ast$, and 
\item by inspection is also present in the coequaliser for $\epsilon = \emptyset$ (since $c' \neq \ast$, so that $\mathrm{F}_\emptyset(c') = \mathrm{F}_\ast(c')$).
\end{enumerate}
This shows that the map \eqref{eq:middleofaproof} is surjective, and injectivity is only notationally more involved, so we leave it to the reader.
\end{proof}

\subsection{The proof of Theorem \ref{thm:duggerspivak}} 
We now have a concrete description of mapping spaces of the kind $\F_{\CC(S)}(s,t)$ in terms of necklaces. The goal of this section is to prove Theorem \ref{thm:duggerspivak}, the statement of which we repeat here:

\begin{theorem}\label{thm:duggerspivakrestate}
Let $S$ be an $\infty$-category and $s,t \in S_0$. Then the adjoint
\[|\Hom^\mathrm R_S(s,t)|_Q \longrightarrow \F_{\CC(S)}(s,t)\]
of the map constructed in Section \ref{subsec:coherentnerve} is a weak homotopy equivalence. 
\end{theorem}

Let us unravel the description of the map
\[|\Hom^\mathrm R_S(s,t)|_Q \longrightarrow \F_{\CC(S)}(s,t)\]
of Theorem \ref{thm:duggerspivakrestate} in terms of necklaces. (Note that we have added the superscript $R$ everywhere to stress that we are working with \emph{right} mapping spaces.) Write $I^n$ for the simplicial set $(\Delta^n \star \Delta^0)/\Delta^n$. Writing $0$ for its initial vertex (the image of $0, \ldots, n$ under the quotient map) and $1$ for its final vertex, the simplicial set $Q^n$ was defined to be $\F_{\CC(I^n)}(0,1)$. By Proposition \ref{prop:Cnecklace} we have
\begin{equation*}
\F_{\CC(I^n)}(0,1) \cong \colim_{\NL \in \Nec / I^n_{0,1}}  \F_{\CC(\NL)}(\alpha(\NL),\omega(\NL)),
\end{equation*}
where $\Nec / I^n_{0,1}$ denotes the full subcategory of $\Nec / I^n$ on necklaces in $I^n$ starting at $0$ and ending in $1$.

Recall that, any simplicial set $X$ can be written as a colimit over its category of simplices as follows:
\[\colim_{\Delta^n \rightarrow X} \Delta^n \cong X.\]
Applying this to $\Hom^\mathrm R_S(s,t)$ we find
\begin{align*}
|\Hom^\mathrm R_S(s,t)|_Q & \cong \colim_{\Delta^n \rightarrow \Hom^\mathrm R_S(s,t)} Q^n \\
& \cong \colim_{I^n \rightarrow S} Q^n \\
& \cong \colim_{I^n \rightarrow S} \colim_{\NL \in \Nec / I^n_{0,1}}  \F_{\CC(\NL)}(\alpha(\NL),\omega(\NL)).
\end{align*}
On the second line, the colimit is over the diagram with objects the maps $I^n \rightarrow S$ sending $0$ to $s$ and $1$ to $t$, with morphisms between them dictated by the cosimplicial structure of $I^\bullet$. The final expression may also be written as a single colimit over the category we will denote $\Nec/I^\bullet/S_{s,t}$. Its objects are of the form $\NL \rightarrow I^n \rightarrow S$, with the first map giving a necklace from $0$ to $1$ in $I^n$ and the second map sending $0$ to $s$ and $1$ to $t$. Clearly there is a forgetful functor
\[\Nec/I^\bullet/S_{s,t} \longrightarrow \Nec/S_{s,t}.\]
This functor induces the top horizontal arrow in the following square, which is easily verified to commute:
\[
\begin{tikzcd}
\colim_{\NL \rightarrow I^n \rightarrow S}  \F_{\CC(\NL)}(\alpha(\NL),\omega(\NL)) \ar{r}\ar{d}{\cong} & \colim_{\NL \rightarrow S}  \F_{\CC(\NL)}(\alpha(\NL),\omega(\NL)) \ar{d}{\cong} \\
\vert\Hom^\mathrm R_S(s,t)\vert_Q \ar{r} & \F_{\CC(S)}(s,t).
\end{tikzcd}
\]

Our first step will be to replace the colimits by homotopy colimits:

\begin{proposition}
\label{prop:colimvshocolim}
The vertical maps in the following squares are weak homotopy equivalences of simplicial sets:
\[
\begin{tikzcd}
N(\Nec/I^\bullet/S_{s,t}) \ar{r} & N(\Nec/S_{s,t}) \\
\hocolim_{\NL \rightarrow I^n \rightarrow S}  \F_{\CC(\NL)}(\alpha(\NL),\omega(\NL)) \ar{r}\ar{d}\ar{u} & \hocolim_{\NL \rightarrow S}  \F_{\CC(\NL)}(\alpha(\NL),\omega(\NL)) \ar{d}\ar{u} \\
\colim_{\NL \rightarrow I^n \rightarrow S}  \F_{\CC(\NL)}(\alpha(\NL),\omega(\NL)) \ar{r} & \colim_{\NL \rightarrow S}  \F_{\CC(\NL)}(\alpha(\NL),\omega(\NL)).
\end{tikzcd}
\]
\end{proposition}
\begin{remark}
For the sake of concreteness, the reader should keep the explicit definition of homotopy colimits via the bar construction in mind, as we will be using it occasionally. For a simplicial diagram $F\colon \C \rightarrow \sSet$ on a category $\C$, the resulting formula is that of Bousfield--Kan, describing the homotopy colimit as the diagonal/realisation of the following bisimplicial set: 
\[(\hocolim_\C F)_{m,n} = \coprod_{c_0 \rightarrow \cdots \rightarrow c_n} F(c_0)_m.\]
The coproduct is over $n$-simplices of the nerve $\nerve \C$, which in particular gives the upper vertical maps in the previous proposition.
\end{remark}
\begin{proof}[Proof of Proposition \ref{prop:colimvshocolim}]
For the upper square this follows from the fact that the simplicial sets $\F_{\CC(\NL)}(\alpha(\NL),\omega(\NL))$ are weakly contractible, so let us focus on the lower square.

The right vertical map in the lower square is a weak equivalence by \cite[Theorem 5.3]{duggerspivakrigid}; in terms of the bisimplicial set of the remark, one finds the $n$th level of the homotopy colimit given by the nerve of the category of necklaces in $S$ with a flag of length $n$ and this is equivalent to its set of components by \cite[Proposition 4.10]{duggerspivakrigid}, which assemble to the actual colimit. To see the claim for the left vertical map in the lower square, we observe that this arrow 
is weakly equivalent to the left vertical arrow in the following commutative square:
\[
\begin{tikzcd}
  \hocolim_{I^n \rightarrow S} \hocolim_{\NL \rightarrow I^n}  \F_{\CC(\NL)}(\alpha(\NL),\omega(\NL)) \arrow[d] \arrow[r, ]  & \hocolim_{I^n \rightarrow S} Q^n\arrow[d] &  \\
  \colim_{I^n \rightarrow S} \colim_{\NL \rightarrow I^n}  \F_{\CC(\NL)}(\alpha(\NL),\omega(\NL))  \arrow[r] & \colim_{I^n \rightarrow S} Q^n. 
\end{tikzcd}
\]
The bottom map in this square is an isomorphism on account of Proposition \ref{prop:Cnecklace}, this time applied to the simplicial set $I^n$ instead of $S$. The top map is a weak homotopy equivalence by another application of \cite[Theorem 5.3]{duggerspivakrigid}. The right map is obtained by applying $|\cdot|_Q$ to the map
\begin{equation*}
\hocolim_{\Delta^n \rightarrow \Hom_S^\mathrm R(s,t)} \Delta^n \rightarrow \colim_{\Delta^n \rightarrow \Hom_S^\mathrm R(s,t)} \Delta^n = \Hom_S^\mathrm R(s,t).
\end{equation*}
But this map is a weak homotopy equivalence for any simplicial set $X$ in place of $\Hom_S^\mathrm R(s,t)$, and $|-|_Q$ preserves weak homotopy equivalences since $Q$ is Reedy cofibrant (see the discussion before Proposition \ref{cor:fiberUn}). 
\end{proof}
\begin{remark}
\label{rmk:necklacehocolims}
Thus, Propositions \ref{prop:colimvshocolim} and \ref{prop:Cnecklace} imply that for a simplicial set $S$ and vertices $s,t \in S_0$, the space $\F_{\CC(S)}(s,t)$ is weakly equivalent to the nerve of the category $\Nec/S_{s,t}$.
\end{remark}

To prove Theorem \ref{thm:duggerspivakrestate} it will now suffice to prove that the top horizontal map of Proposition \ref{prop:colimvshocolim} is a weak homotopy equivalence. 

\begin{proposition}
\label{prop:weqnerves}
Let $S$ be an $\infty$-category and $s,t \in S_0$. Then the functor
\[\Nec/I^\bullet/S_{s,t} \longrightarrow \Nec/S_{s,t}, \quad  (\NL \rightarrow I^n \rightarrow S) \longmapsto (\NL \rightarrow S)\]
induces a weak homotopy equivalence between the nerves of these categories.
\end{proposition}

The remainder of this section is devoted to proving this result. We introduce an auxiliary category (related to what Dugger and Spivak call \emph{gadgets}):

\begin{definition}
A \emph{categorical interval} is a doubly pointed simplicial set $X$ such that both mapping complexes
\[\F_{\CC(X)}(0, 1) \text{ and } \Hom_{X_f}^\mathrm{R}(0, 1)\]
are weakly contractible, where $X_f$ is any fibrant replacement of $X$ in the Joyal model structure on $\sSet_{*,*}$. We write $\mathcal{G}$ for the full subcategory of $\sSet_{\ast, \ast}$ spanned by the categorical intervals.
\end{definition}

This category serves as a convenient collection of simplicial sets that contains both necklaces and the $I^n$:

\begin{lemma}
Each $I^n$ and every necklace is a categorical interval.
\end{lemma}
\begin{proof}
For necklaces one observes that the inclusion
\[
\Delta^{n_0} \vee \dotsb \vee \Delta^{n_k} \rightarrow \Delta^{n_0 + \cdots + n_k}
\]
is inner anodyne. It is clear that the codomain is a categorical interval; since the functor $\CC$ sends categorical equivalences to equivalences of simplicial categories (Proposition \ref{cor:CCQuillen}), it follows that the domain is a categorical interval as well. For $I^n$, consider the pushout square
\[
\begin{tikzcd}
\Delta^n \vee \Delta^1 \ar{d}\ar{r} & \Delta^n \star \Delta^0 \ar{d} \\
\Delta^1 \ar{r} & I^n.
\end{tikzcd}
\]
The top horizontal arrow is inner anodyne, hence so is the bottom. Since $\Delta^1$ is a categorical interval, the same is true of $I^n$.
\end{proof}

We write $\mathcal{G}/S_{s,t}$ for the slice category of categorical intervals over $S$ sending $0$ to $s$ and $1$ to $t$. The inclusions of $I^\bullet$ and $\Nec$ into $\mathcal{G}$ now give rise to a square of functors
\[
\begin{tikzcd}
\Nec/I^\bullet/S_{s,t} \ar{r}\ar{d} & \Nec/S_{s,t} \ar{d}{\mathrm{incl}} \\
I^\bullet/S_{s,t} \ar{r}{\mathrm{incl}} & \mathcal{G}/S_{s,t}.
\end{tikzcd}
\]
The two unlabelled arrows are the forgetful functors. Note that this square does \emph{not} commute; however, there is an evident natural transformation from the composite along the right towards the composite along the left. Thus, the square commutes up to homotopy after passing to classifying spaces. Moreover, it is not hard to see that the vertical arrows induce weak homotopy equivalences:

\begin{lemma}
\label{lem:inclsweqs}
The inclusion $\Nec/S_{s,t} \rightarrow \mathcal{G}/S_{s,t}$ and the forgetful functor $\Nec/I^\bullet/S_{s,t} \rightarrow I^\bullet/S_{s,t}$ induce weak homotopy equivalences of classifying spaces.
\end{lemma}
\begin{proof}
For a fixed object $f\colon J \rightarrow S_{s,t}$ of the category $\mathcal{G}/S_{s,t}$, the slice of the forgetful functor $\Nec/S_{s,t} \rightarrow \mathcal{G}/S_{s,t}$ over $f$ can be identified with the category $\Nec/J_{0,1}$. By Remark \ref{rmk:necklacehocolims}, the classifiying space of this category is weakly equivalent to $\CC(J)(0,1)$, which is weakly contractible because $J$ is a categorical interval. Hence Quillen's Theorem A implies and the functor $\Nec/S_{s,t} \rightarrow \mathcal{G}/S_{s,t}$ induces a weak homotopy equivalence. The same argument applies to the other functor $\Nec/I^\bullet/S_{s,t} \rightarrow I^\bullet/S_{s,t}$, now relying on the fact that $\CC(I^n)(0,1)$ is weakly contractible for any $n$.
\end{proof}

From Lemma \ref{lem:inclsweqs} and the homotopy commutative square above it, we see that the proof of Proposition \ref{prop:weqnerves} is completed by the following result:

\begin{proposition}
\label{prop:catintervalsclassifyingspace}
The inclusion $I^\bullet/S_{s,t} \rightarrow \mathcal{G}/S_{s,t}$ induces a weak homotopy equivalence of classifying spaces.
\end{proposition}

Note that we cannot directly apply Quillen's Theorem A to show this. Indeed, the relevant slice categories of the form $I^\bullet/X_{0,1}$, for some categorical interval $X$, would be of the homotopy type of the right mapping space $\Hom_{X}^\mathrm R(0, 1)$. This can only be expected to be weakly contractible in case $J$ is an $\infty$-category. Therefore, let us introduce the full subcategory $\mathcal{G}_f$ of $\mathcal{G}$ spanned by the fibrant categorical intervals. We will exploit the following observation:

\begin{lemma}
The inclusion $\mathcal{G}_f/S_{s,t} \subset \mathcal{G}/S_{s,t}$ induces a weak homotopy equivalence of classifying spaces.
\end{lemma}
\begin{proof}
Let us denote this inclusion of categories by $i$. Pick a fibrant replacement functor $$R\colon \sSet_{\ast, \ast}/S_{s,t} \rightarrow \sSet_{\ast, \ast}/S_{s,t}$$ with respect to the Joyal model structure,
which comes with a natural transformation $\eta\colon \mathrm{id} \Rightarrow R$. Since $\mathcal{G}$ is closed 
under categorical equivalences and $S$ is an $\infty$-category, the functor $R$ restricts to a functor $\mathcal{G}/S_{s,t} \rightarrow \mathcal{G}_f/S_{s,t}$, and the natural transformation 
$\eta$ restricts to natural transformations $\mathrm{id}_{\mathcal{G}_f/S_{s,t}} \Rightarrow Ri$ and $\mathrm{id}_{\mathcal{G}/S_{s,t}} \Rightarrow iR$. Upon passing to classifying spaces, these natural 
transformations witness 
$i$ as a homotopy equivalence with $R$ as homotopy inverse. 
\end{proof} 

We can now finish our argument by the following:

\begin{proof}[Proof of Proposition \ref{prop:catintervalsclassifyingspace}]
We may factor the canonical map $I^\bullet \rightarrow \Delta^1$ of cosimplicial objects 
as a trivial Reedy cofibration $I^\bullet \rightarrow I^\bullet_f$ followed by a (necessarily trivial) Reedy fibration 
$I^\bullet_f \rightarrow \Delta^1$ (with the Reedy model structure based on the sliced Joyal model structure on $\sSet_{\ast,\ast}$). In particular, each term $I^n_f$ is a fibrant categorical interval. Furthermore, we note that the induced map 
$$
\Hom_{\sSet_{\ast, \ast}}(I^\bullet_f, S) \rightarrow \Hom_{\sSet_{\ast, \ast}}(I^\bullet, S)
$$
is a weak homotopy equivalence for every fibrant doubly pointed simplicial set $S$ (see \cite[Proposition 17.1.6]{Hirschhorn} for instance; both sides can be interpreted as the `mapping space' from $\Delta^1$ to $S$ with respect to the Joyal model structure). Consequently, the induced map between the respective categories of simplices $I^\bullet_f/S_{s,t}$ and $I^\bullet/S_{s,t}$ is a weak homotopy equivalence. (Indeed, any simplicial set has the same weak homotopy type as its category of simplices; the standard argument is a skeletal induction, reducing to the case of a simplex, where both are weakly contractible.) The map of cosimplicial objects $I^\bullet \rightarrow I^\bullet_f$ gives rise to a natural transformation
\[
  \begin{tikzcd}
    I^\bullet_f/S_{s,t} \arrow[r] \arrow[d] & \mathcal{G}_f/S_{s,t} \arrow[d] \\ 
    I^\bullet/S_{s,t}  \ar[ur,shorten <=10pt, shorten >=10pt, Rightarrow]  \arrow[r] & \mathcal{G}/S_{s,t},
  \end{tikzcd}
\]
where all arrows are inclusions. We have already demonstrated that the vertical arrows are weak homotopy equivalences (see the previous lemma for the arrow on the right). Consequently, it suffices to show that the inclusion $I^\bullet_f/S_{s,t} \subset \mathcal{G}_f/S_{s,t}$ 
is a weak homotopy equivalence. Consider any object $f\colon X \rightarrow S_{s,t}$ of $\mathcal{G}_f/S_{s,t}$. Then the slice category of the inclusion $I^\bullet_f/S_{s,t} \subset \mathcal{G}_f/S_{s,t}$ over $f$ is precisely $I^\bullet_f/X_{0,1}$. As above, this category is weakly equivalent to the category $I^\bullet/X_{0,1}$, which in turn is the category of simplices of the simplicial set $\Hom_{X}^\mathrm{R}(0, 1)$. The latter is contractible, since $X$ is a fibrant categorical interval. Quillen's Theorem A now implies that the inclusion $I^\bullet_f/S_{s,t} \subset \mathcal{G}_f/S_{s,t}$ is indeed a weak homotopy equivalence.
\end{proof}

\bibliographystyle{amsalpha}


\end{document}